\documentclass{notices}
\usepackage[utf8]{inputenc}
\usepackage{amssymb,amsmath,amsthm,color,mathtools}

\usepackage{tikz,tkz-euclide}
\usepackage[breakable,skins]{tcolorbox}
\newtcolorbox{tbox}[1][]{%
    breakable,
    enhanced,
    colframe=blue,
    coltitle=white,
    #1
}

\hypersetup{
    colorlinks=true,       % false: boxed links; true: colored links
    linkcolor=blue,          % color of internal links
    citecolor=blue,        % color of links to bibliography
    filecolor=blue,      % color of file links
    urlcolor=blue           % color of external links
}

\usepackage{tcolorbox}
\usepackage[all]{xy}

\newtheorem{theorem}{Theorem}[section]
\newtheorem{proposition}[theorem]{Proposition}

\theoremstyle{definition}

\newtheorem{definition}[theorem]{Definition}
\newtheorem{remark}[theorem]{Remark}
\newtheorem{example}[theorem]{Example}

\renewcommand{\div}{\operatorname{div}}
\newcommand{\Eff}{\operatorname{Eff}}
\newcommand{\Mov}{\operatorname{Mov}}
\newcommand{\Nef}{\operatorname{Nef}}
\newcommand{\Cl}{\operatorname{Cl}}

\newcommand{\Bl}{\operatorname{Bl}}

\title{Finite generation of Cox rings}

\date{}
\author{
   Jos\'e Luis Gonz\'alez
  \affil{University of California, Riverside. Email: 
  joselg@ucr.edu
    }
    \thanks{J. Gonzalez was supported by a grant from the Simons Foundation (Award Number 710443).}
  \and
  Antonio Laface
  \affil{
  Universidad de Concepci\'on.
   Email: alaface@udec.cl
   }
   \thanks{A. Laface has been partially supported by Proyecto FONDECYT Regular
n. 1190777.}
}

\usepackage{xcolor}

\begin{document}

\maketitle

\setcounter{tocdepth}{2}

The projective plane $\mathbb P^2$ is one of the most ubiquitous objects  in geometry.
It appears as a compactification 
of the plane by adding 
``a line at infinity'' as proposed 
around the second half of the 15th century
by Piero della Francesca in his 
\href{http://bibdig.museogalileo.it/Teca/Viewer?an=1040190&pag=Carta: 8r}
{\em De prospectiva pingendi}. 
Painters started using perspective as a consequence of this construction.

\begin{figure}[ht]
\begin{center}
\includegraphics[scale=.3]{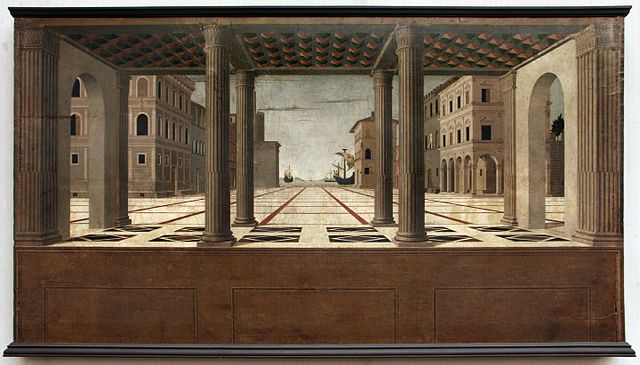}
\caption{Architectural Veduta, 
by Francesco di Giorgio Martini,
around 1490.}
\label{fig1}
\end{center}
\end{figure}

\vspace{-3mm}

One has to wait until the beginning 
of the XIX century for August 
Ferdinand Möbius to introduce 
{\em homogeneous coordinates}
which provide an algebraic framework
for %(Desarguesian)
projective planes.
Homogeneous coordinates appear
as soon as one defines the $n$-dimensional complex projective space as
\[
 \mathbb P^n
 :=
 \left(\mathbb C^{n+1}\setminus\{0\}\right)/\sim 
\] 
where $\sim$ is the equivalence relation that identifies $x=(x_0,\ldots,x_n)$ and $\lambda x = (\lambda x_0,\ldots,\lambda x_n)$ for all $x \in\mathbb C^{n+1}$ and $\lambda \in\mathbb C^* := \mathbb C\setminus\{0\}$. 
{\em Cox~rings} can be 
thought of as a generalization of 
homogeneous coordinates to a 
wider class of algebraic varieties,
like the cartesian products of
projective spaces. 
The main idea behind 
Cox rings is that they allow us to define the
quotient construction of the projective 
space in an intrinsic way: where is the
space $\mathbb C^{n+1}$ coming from?
where is the equivalence relation coming from?
why does one remove the zero vector? 
We provide answers to these questions. 
The projective space $\mathbb P^n$ is the 
disjoint union of the affine space $\mathbb C^n$
with the ``hyperplane at infinity'' $H$, 
given by the horizon line in  
Figure~\ref{fig1}.
The coordinate ring of $\mathbb C^n$ is the 
polynomial ring $S := \mathbb C[x_1,\dots,x_n]$
in $n$ variables.
Given an integer $d\geq 0$ one can form the
vector space $S_d\subseteq S$ of polynomials of
degree at most $d$. 
The Cox ring of $\mathbb P^n$ is the 
direct sum $\oplus_{d\geq 0}S_d$ which is a graded 
algebra isomorphic to the polynomial
ring $\mathbb C[x_0,\dots,x_n]$ in $n+1$ 
variables. This is the coordinate ring of
$\mathbb C^{n+1}$ and its grading assigning degree one to each variable corresponds to the
$\mathbb C^*$-action $\lambda\cdot x := \lambda x$, 
which induces the previous equivalence relation.
Finally, the point $0\in\mathbb C^{n+1}$ is in the
closure of any orbit for this action, then one
has to remove it in order to obtain a quotient which
is not just one point.

From the construction of the Cox ring
of the projective space, we can see that the degree of 
a homogeneous polynomial is a central notion. This degree can be replaced with the order of the pole of a polynomial in $S$ along the ``hyperplane at infinity'' $H$. In both cases, the notion depends on the choice of an open subset of the projective space: the affine space $\mathbb C^n$.
The reason for this choice is that the complement
$H = \mathbb P^n\setminus\mathbb C^n$ generates the
{\em divisor class group} of $\mathbb P^n$, an invariant
which is hidden behind the whole construction
of the Cox ring.
So, more generally, when one starts with an algebraic
variety $X$ which admits an affine  Zariski open subset
$U\subseteq X$ such that the 
codimension-one subvarieties 
in $X\setminus U$ generate the divisor class group of
$X$ one can define the Cox ring of $X$ in a similar
way as for the Cox ring of $\mathbb P^n$.
The Cox ring is not necessarily finitely generated as an algebra over the complex numbers, but 
when it is finitely generated, it is the coordinate ring of
an affine algebraic variety $\overline X$.
The grading of the Cox ring corresponds to an action 
of an abelian group $H_X$ over $\overline X$.
One can show that, as in the case of the projective space,
there is an open invariant subset $\widehat X$ 
of $\overline X$ and a quotient map
\[
 p_X\colon \widehat X\to X
\]
by the action of $H_X$.

Understanding all morphisms from a projective variety $X$ to other projective varieties is a fundamental problem. One would like to decompose each map into simple steps and parametrize the possibilities. 
Morphisms between projective varieties   
can be factored as the composition of a morphism with connected fibers followed by one with finite fibers. 
Then, we would like to understand the decompositions of maps with connected fibers. 
A natural setting for this question is to consider not only morphisms but also the so-called {\em rational contractions}, that is,  
compositions of rational maps which are isomorphisms in codimension one and surjective morphisms with
connected fibers, see \cite{HK}*{Definition 1.0}. 
If the Cox ring of $X$ is finitely generated 
each such rational contraction
is uniquely determined by the choice of  
a certain invariant open subset of $\overline X$, 
which gives rise to another quotient.
The collection of these subsets
is in bijection with the set of cones,
called {\em Mori chambers},
of a fan supported on the \emph{effective cone} of $X$,  
see Section~\ref{section.mori.chambers}.

The question of the finite generation
of the Cox ring is thus central and
we  devote the second part of this
note to it, ending it with some recent 
examples of finitely and non-finitely 
generated Cox rings. 
Along the way, we  discuss some
natural questions, including when are
Cox rings polynomial rings?
when do two algebraic varieties 
have isomorphic Cox rings?
The answer to the first question
leads us to \emph{toric varieties} and
the work of D. Cox~\cite{Cox},
while the answer to the second 
question allows us to introduce
a fascinating geometric object:
{\em Mori dream spaces}, following
the work of Y. Hu and S. Keel~\cite{HK}.

A comprehensive reference on Cox rings is~\cite{ADHL}, which also discusses the history of the subject. We highlight J.-L. Colliot-Th\`{e}l\'ene and J.-J. Sansuc who introduced universal torsors in arithmetic geometry in the 1970s. David Cox introduced the homogeneous coordinate rings of toric varieties in~\cite{Cox}. Y. Hu and S. Keel proposed the name Cox ring and showed how the finite generation of this ring is connected to the geometry of the variety. Finally, J. Hausen in~\cite{hausen2008cox} gave a definition
of Cox rings which generalizes the previous ones
and is the one described in this article.

%%%%%%%%%%%%%%%%%%%%%%%%%%%%%%%%%%%%%%%%%%%%%%%%%%

\section{Graded algebras}   \label{section.graded.algebras}
Our first step towards the definition
of Cox rings is to introduce
graded algebras.
Given a finitely generated abelian 
group $A$, 
a $\mathbb{C}$-algebra $R$ is {\em $A$-graded} if 
\[
 R := \bigoplus_{a\in A}R_a,
\]
where each $R_a$ is a complex vector
space and $R_{a}\cdot R_{b}\subseteq 
R_{a+b}$ for any $a,b\in A$.
The elements of $R_a$ are called 
{\em homogeneous} of degree $a$.
One says that $R$ is {\em finitely 
generated}
if there is a surjective homomorphism 
of $\mathbb C$-algebras from the polynomial ring 
$\mathbb C[x_1,\dots,x_r]$ to $R$.
Moreover, without loss of generality,
one can assume the images of the 
variables are homogeneous.
The kernel is the ideal of an affine algebraic set $X \subseteq \mathbb C^r$ which is uniquely determined by $R$ up to isomorphism.
One can recover $R$ as the coordinate ring of $X$, and this ring is a domain precisely when $X$ is an affine variety (i.e., an irreducible affine algebraic set).
The group $A$ determines the {\em monoid algebra}
\[
 \mathbb C[A]
 :=
 \bigoplus_{a\in A}\mathbb C\cdot \chi^a
\]
with product $\chi^a\cdot \chi^b = \chi^{a+b}$  (see~\cite{ADHL}*{Constr. 1.1.1.5}).
One can show that the affine variety
$G$ defined by this algebra is 
isomorphic to ${\rm Hom}(A,\mathbb C^*)$
as a group. 
If $A$ is isomorphic to 
$\mathbb{Z}^n\oplus A_{\rm tor}$, where the
second summand denotes the torsion
part, then we have 
$G\simeq (\mathbb C^*)^n\oplus A_{\rm tor}$.
The latter group is a {\em quasitorus}
which acts on $X$ in the following way.
For any $a\in A$ denote by $\chi^a
\colon G\to\mathbb C^*$ the character
$g\mapsto \chi^a(g) := g(a)$.
If we denote by $a_1,\dots,a_r\in A$
the degrees of the above $r$ homogeneous
generators of $R$, then there is an 
action of $G$ on $\mathbb C^r$ given by 
\[
 g\cdot (x_1,\dots,x_r)
 := 
 (\chi^{a_1}(g)x_1,\dots,\chi^{a_r}(g)x_r), 
\]
which induces an action of $G$ on $X$.
In more abstract terms, this group action
is induced by the homomorphism $R\to 
R\otimes\mathbb C[A]$ defined by 
$R_a\ni f_a\mapsto f_a\otimes\chi^a$.
The above action can
be defined even when $R$ is not finitely
generated.

\begin{example}
\label{ex-1}
Two examples, useful in the sequel, 
are given by two different $\mathbb Z^2$-gradings 
of the polynomial ring $\mathbb C[x_1,x_2,x_3,x_4]$.
Each grading is defined by the columns of one 
of the following matrices 
(the degree of $x_i$ is the 
$i$-th column):
\[
 \begin{bmatrix}
  1&1&n&0\\
  0&0&1&1
 \end{bmatrix}
 \qquad
 \qquad
 \begin{bmatrix}
  1&1&1&0\\
  0&1&1&1
 \end{bmatrix},
\]
where $n$ is a nonnegative integer.
In each case $G \simeq(\mathbb C^*)^2$ acts on $\mathbb C^4$. 
If we denote by
$t_1 := \chi^{(1,0)}(t)$ and 
$t_2 := \chi^{(0,1)}(t)$ then, for example, 
the second action is 
$t\cdot (x_1,\dots,x_4) = 
(t_1x_1,t_1t_2x_2,t_1t_2x_3,t_2x_4)$.
\end{example}

%%%%%%%%%%%%%%%%%%%%%%%%%%%%%%%%%%%%%%%%%%%%%%%%%%

\subsubsection*{Divisor class group}
Our second step is to introduce a
fundamental invariant of an algebraic
variety: the {\em divisor class group}.
Let $X$ be an algebraic variety.
A rational map or rational function $X\dashrightarrow \mathbb C$ on $X$ is an equivalence class of morphisms from nonempty open subsets of $X$ to $\mathbb{C}$, where two such morphisms are identified if they agree on a nonempty open set.
The set of all rational functions on $X$ is a field under pointwise sum and product 
called the {\em field of 
rational functions} of $X$. 
The free abelian group generated 
by~the codimension-one subvarieties
is the {\em group of Weil divisors} of $X$. 
Given a Weil divisor $D = \sum_ia_iD_i$,
its {\em support} is the union $\bigcup_iD_i$. 
A codimension-one sub\-va\-ri\-e\-ty is called
a {\em prime divisor.}
To any nonzero rational function
$f\in\mathbb C(X)$ one associates
a Weil divisor 
\[
 \div(f) := \sum_{D\text{ prime divisor}}{\rm ord}_D(f)D,
\]
where ${\rm ord}_D(f)$ is an integer which
intuitively represents the order of vanishing
of $f$ along the prime divisor
$D\subseteq X$. It turns out that the sum above is finite and that the assignment 
$f\mapsto \div(f)$ is a homomorphism 
of abelian groups. Its image is the group
of {\em principal divisors} and finally one 
defines the divisor class group ${\rm Cl}(X)$
as the quotient of the group of Weil divisors 
modulo the subgroup of principal divisors.
From now on we  assume that this group
is finitely generated, which 
is the case for example
when $X$ is a rational, 
a Fano or a Calabi-Yau
variety.
Two Weil divisors $D, D'$ are 
{\em linearly equivalent} if their difference $D-D'$ is principal.
We now recall the definition of a normal variety.

\begin{tbox}
\begin{definition}
A {\em normal affine variety} $X$
is an affine variety which has
singularities in codimension 
two or more and that for any open subset 
$U\subseteq X$, with complement of codimension
two or more, the regular functions of 
$U$ are restrictions of regular 
functions of $X$.
More generally a normal variety is one
which is covered by normal affine ones.
\end{definition}
\begin{example}
An example of non-normal affine variety 
is given by any singular curve, since
the first condition for normality is not satisfied.
\end{example}
\begin{example}
An example where the second condition is
not satisfied is the affine variety $X$ defined
by the subalgebra $A\subseteq\mathbb C[x,y]$
generated by all the monomials but $x$.
It is easy to see that $A$ is in fact 
generated by $x^2,x^3,y,xy$, so that $X$ 
is the image in $\mathbb C^4$ of the map
$(x,y)\mapsto (x^2,x^3,y,xy)$. 
If $u_1,u_2,u_3,u_4$ are coordinates of
$\mathbb C^4$ then on the open subset
$U := X\setminus V(u_1,u_3)$ the function
$u_2/u_1 = u_4/u_3$ is regular but it is 
not restriction of a regular function on $X$.
\end{example}
\end{tbox}

Our next definition is that of a {\em Cartier
divisor}. The usual definition is different
from the one given in the following lines,
but one can show that on a normal algebraic
variety the two definitions coincide.
A Weil divisor $D$ is {\em Cartier} if 
it is locally principal, that is 
$X$ admits an open covering and 
over each such open subset $U$ 
one has $D|_U = \div(f)|_U$ for
some rational function $f\in\mathbb C(X)$,
where $D|_U$ means that one removes all the
prime divisors in the support of $D$ that
do not intersect $U$.
A variety is {\em $\mathbb Q$-factorial} 
if any Weil divisor admits a nonzero
integer multiple which is a Cartier
divisor.

%%%%%%%%%%%%%%%%%%%%%%%%%%%%%%%%%%%%%%%%%%%%%%%%%%

%\subsubsection*{Sheaves}
As a preliminary to the definition of Cox rings we briefly recall the notion of a sheaf.

\begin{tbox}
\begin{definition}
A {\em presheaf} $\mathcal{F}$
of groups over a topological space $X$ 
is a contravariant functor from the 
category of open subsets of $X$ with inclusions
to the category of groups  
with homomorphisms.
To any inclusion of open sets 
$V\subseteq U$ one associates a homomorphism
$\mathcal F(U)\to \mathcal F(V)$ usually denoted by
$f\mapsto f|_V$.
The set $\mathcal{F}(U)$ is also denoted by $\Gamma(\mathcal{F},U)$ and its elements are called the sections of $\mathcal{F}$ over $U$. 
A presheaf is a {\em sheaf} if for any 
open covering $\{U_i\}$ of an open
subset $U$ and any collection of 
$f_i\in\mathcal F(U_i)$, such that 
${f_i}|_{U_i\cap U_j} = {f_j}|_{U_i\cap U_j}$
for all $i,j$, there exists a unique 
$f\in\mathcal{F}(U)$ such that 
$f|_{U_i} = f_i$ for all $i$.
Presheaves and sheaves of rings or algebras are defined analogously. 

\begin{example}
A typical example to keep in mind
is the presheaf of constant functions
with real values on a topological space.
It is not a sheaf in general since one can have
locally constant functions over disjoint
connected components which do not form
a constant function on the whole space
(here the sheaf condition is automatic
being the intersections $U_i\cap U_j$
all empty).
The presheaf of locally constant functions 
is a sheaf.
\end{example}

\begin{example}
Two important examples are the sheaves that assign to each open subset $U$ of a variety $X$ the regular functions on $U$ and the rational functions on $U$, with restrictions given by the restriction of functions. 
\end{example}

\end{definition}
\end{tbox}

%%%%%%%%%%%%%%%%%%%%%%%%%%%%%%%%%%%%%%%%%%%%%%%%%%

A Weil divisor $D = \sum_ia_iD_i$ 
is {\em effective} if all of its coefficients $a_i$
are nonnegative. Any Weil divisor $D$ defines
a sheaf $\mathcal O_X(D)$ on $X$ whose 
space of sections on the open subset $U\subseteq X$
is the vector space generated by the nonzero 
rational functions $f\in\mathbb C(X)$ such
that $(\div(f)+D)|_U$ is effective. 
A basis $\{f_0,\dots,f_N\}$ of the space
of global sections $\Gamma(X,\mathcal O_X(D))$
induces a rational map
\[
 \psi_D\colon 
 X\dashrightarrow\mathbb P^N,
 \qquad
 p\mapsto [f_0(p):\cdots:f_N(p)],
\]
and choosing a different basis changes 
$\psi_D$ with $\sigma\circ\psi_D$,
where $\sigma$ is an automorphism of 
$\mathbb P^N$.
The divisor $D$ is {\em base point free} if $\psi_{D}$ is a morphism;
it is {\em very ample} if
$\psi_D$ is an embedding; it is 
{\em ample}, respectively 
{\em semiample}, if there exists 
positive integer $n$ such that 
$\psi_{nD}$ is very ample, 
respectively $\psi_{nD}$ is a morphism.

%%%%%%%%%%%%%%%%%%%%%%%%%%%%%%%%%%%%%%%%%%%%%%%%%%

\subsubsection*{Cox rings}
To any finitely generated subgroup $K$ of Weil divisors on a normal algebraic variety $X$ one can associate the following sheaf of algebras
\begin{equation}
\label{defS}
 \mathcal S := \bigoplus_{D\in K}\mathcal O_X(D).
\end{equation}
The sheaf $\mathcal S$ is graded
by $K$ and it is possible to show 
that its ring of global sections
$\Gamma(X,\mathcal S)$ is a 
factorial ring
whenever the class map 
${\rm cl}|_K\colon K\to 
{\rm Cl}(X)$ is a surjection
~\cite{ADHL}*{Thm. 1.3.3.3}.
If the divisor class group is torsion-free
then one can choose $K$ such that
${\rm cl}|_K$ is an isomorphism 
and define a Cox ring of $X$ as
$\Gamma(X,\mathcal S)$. 
It is not difficult to show that different
choices for $K$ lead to isomorphic 
Cox rings. 
However things get more complicated
when the divisor class group has torsion. 
In order to provide a definition
of Cox ring which also includes this
possibility one starts with a 
finitely
generated subgroup of Weil divisors 
$K$ such that ${\rm cl}|_K$ is surjective and
let $K_0$ be its kernel. 
One can define a homomorphism
$\chi\colon K^0\to \mathbb C(X)^*$
such that $\div\circ\chi = {\rm id}$
because $K^0$ is free abelian and its elements are principal 
divisors. 
Given a principal divisor $D$, 
a rational function $f$ such that
$\div(f) = D$ is defined only up to 
scalar multiplication, provided
that the only global invertible regular
functions of $X$ are constants, which 
we are going to assume from now on. Introducing 
$\chi$ allows one to make a coherent
choice for all these scalars.
Define the sheaf of ideals
$\mathcal I\subseteq\mathcal S$
which is locally generated by elements 
of the form $1-\chi(D)$, for $D\in K^0$.
Given two divisors $D,D' \in K$ whose 
difference is in $K^0$ it is easy to see
that the map $\mathcal O_X(D)\to
\mathcal O_X(D')$ defined by multiplication
by $\chi(D-D')$ is an isomorphism.
Taking the quotient by $\mathcal I$
has the effect of identifying these 
sheaves keeping one copy for each
divisor class.
\begin{definition}   \label{definition.cox.ring} 
Let $X$ be a normal algebraic variety 
with a finitely generated divisor class
group and whose global invertible regular
functions are constants.
Given a choice of $K$ and $\chi$ as
before one defines a {\em Cox sheaf} 
and a {\em Cox ring} of $X$ as
\[
 \mathcal R := \mathcal S/\mathcal I
 \quad
 \text{and}
 \quad
 \mathcal R(X) := \Gamma(X,\mathcal R),
\]
respectively.
It is possible to show that any two such 
Cox rings for $X$ are isomorphic, in other
words the isomorphism class does not
depend on the choice of $K$ and $\chi$.
\end{definition}

We can write the Cox ring making explicit its grading by the divisor 
class group as follows:
\[
 \mathcal R(X) = \bigoplus_{{\rm Cl}(X)}
 \Gamma(X,\mathcal O_X(D)).
\]

\begin{remark}   \label{remark.small.modification}
A {\em big open subset} of an algebraic variety is an open subset whose complement
has codimension greater than or equal to $2$. 
A direct consequence of Definition~\ref{definition.cox.ring} is that if $U\subseteq X$
is a big open subset, then the 
inclusion induces an isomorphism between their   
Cox rings. 
Indeed, each prime divisor of $X$ 
restricts to a prime divisor of 
$U$ and conversely a divisor of $U$ has
a unique closure in $X$.
In particular, if $X$ and $Y$ are birationally
equivalent varieties with isomorphic big
open subsets then the Cox rings of $X$
and $Y$ are isomorphic.
\end{remark}

\begin{example}
Consider %the projective space 
$X = \mathbb P^n$. The divisor class group
is freely generated by the class of a hyperplane
$H$, so let $K = \langle H\rangle$.
If $x_0,\dots,x_n$ are homogeneous
coordinates and $H = V(x_0)$, 
then the space of global sections of 
$\mathcal O_{\mathbb P^n}(dH)$ is
generated by polynomials in the variables
$\frac{x_1}{x_0}, \dots,\frac{x_n}{x_0}$
of degree up to $d$. This space is isomorphic
to the space of degree $d$ homogeneous
polynomials in $n+1$ variables.
It follows that the Cox
ring is isomorphic to the polynomial
ring $\mathbb C[x_0,\dots,x_n]$.
\end{example}

We conclude this section with an algebraic 
property of Cox rings. 
If the divisor 
class group is torsion-free then the Cox ring
is a factorial ring.
More generally one has the following,
see \cite{ADHL}*{Thm. 1.5.3.7, Thm. 1.3.3.3}
and \cites{berchtold2003homogeneous, ELIZONDO2004625, arzhantsev2009factoriality}.
\begin{theorem}
Every nonzero non-unit homogeneous element in
the Cox ring can be written as a product of a
finite number of irreducible homogeneous
elements, uniquely up to order and units.
\end{theorem}

\section{Quotient construction}
Whenever the Cox ring of an algebraic
variety $X$ 
is finitely generated it defines an affine 
variety named the 
{\em total coordinate space} 
of $X$ and denoted by $\overline X$.
It follows, by our previous discussion, 
that the monoid algebra $\mathbb C[{\rm Cl}(X)]$
defines a quasitorus $H_X$ which
is isomorphic, as a group, to ${\rm Hom}({\rm Cl}(X),\mathbb C^*)$. 
The ${\rm Cl}(X)$-grading of the Cox ring
induces an action of $H_X$ on $\overline X$.
The Cox sheaf determines an invariant open 
subset $\widehat X\subseteq\overline X$
which admits $X$ as good quotient.
We briefly describe how 
this invariant open subset is constructed.
Let $U_1,\dots,U_s$ be affine open subsets
whose union is $X$.
The algebra of global sections 
$\Gamma(U_i,\mathcal R)$ defines 
an invariant affine subvariety 
$\overline U_i\subseteq\overline X$
and one defines $\widehat X := \overline U_1\cup\cdots
\cup\overline U_s$.
The inclusion homomorphism
$\Gamma(U_i,\mathcal R)_0\subseteq
\Gamma(U_i,\mathcal R)$ of the subring 
generated by the degree-zero homogeneous 
elements induces a quotient map
$\overline U_i \to U_i$. All these maps
glue together producing the quotient map
$p_X\colon \widehat X\to X$.
Thus one has the following diagram 
\[
 \xymatrix@C=2pt{
  \widehat X\ar[d]^-{p_X} & \subseteq & \overline X\\
  X.
 }
\]
\begin{example}
When $X$ is the projective space 
$\mathbb P^n$ we have already 
seen that the Cox ring is the polynomial
ring $\mathbb C[x_0,\dots,x_n]$ with the 
$\mathbb Z$-grading which assigns 
degree $1$ to each variable. 
One can cover the projective space with the
$n+1$ affine spaces $U_0,\dots,U_n$, where
$U_i := \mathbb P^n\setminus V(x_i)$.
Over each such subset
$\Gamma(U_i,\mathcal R)
 \simeq
 \mathbb C[x_0,\dots,x_n]_{x_i}$,
where the right-hand side is 
the localization of the polynomial
ring with respect to the multiplicative
subset generated by the powers of $x_i$.
In other words, one takes Laurent 
polynomials in the $i$-th variable.
Thus, for example, $\overline U_0$ 
is the open subset of $\mathbb C^{n+1}$
where the first variable does not vanish.
The quotient morphism 
$\overline U_0\to U_0$ is given by 
$(x_0,\dots,x_n)\mapsto
[1:x_{1}/x_0:\dots:x_n/x_0]$.
All these morphisms glue
together to give the usual quotient 
construction of the projective space.
\end{example}

\begin{remark}
The morphism $p_X$ is a {\em good
quotient} in the following sense:  it is {\em invariant}, which means that 
$p_X(h\cdot \hat x) = p_X(\hat x)$
for any $h\in H_X$ and $\hat x\in\widehat X$; 
it is {\em affine}, which means that the preimage 
of an affine subset $U$ is affine;  
and finally each invariant regular function 
on the affine subset $p_X^{-1}(U)$ is the pullback 
of a regular function on $U$.
\end{remark}

\begin{example}
\label{ex-2}
If we go back to Example~\ref{ex-1}
the two gradings on $R := \mathbb C[x_1,x_2,x_3,x_4]$
define two actions of $(\mathbb C^*)^2$ on
$\mathbb C^4$. 
Let $U_{ij} := \mathbb C^4\setminus V(x_ix_j)$
for $i\in\{1,2\}$ and $j\in\{3,4\}$ 
and let $U$ be the union of these
four open affine subsets.
In both cases the grading 
induces a quotient map of tori 
$(\mathbb C^*)^4\to(\mathbb C^*)^2$.
In the first case it is
$(x_1,x_2,x_3,x_4)\mapsto
(x_2/x_1,x_1^nx_4/x_3)$
and the four quotient maps glue together 
to give the good quotient
$U\to\mathbb F_n$, where 
$\mathbb F_n$ is the $n$-th Hirzebruch 
surface.
In the second case the morphism is
$(x_1,x_2,x_3,x_4)\mapsto
(x_2/x_3,x_1x_4/x_3)$
and the images of $U_{13}$ and 
$U_{34}$ coincide, as a consequence 
the quotient morphism 
on $U$ is not affine.
\end{example}

Example~\ref{ex-2} shows that not every
graded polynomial ring is a Cox ring. 
A complete characterization is given by 
the next proposition, where $A_{\mathbb Q}$ denotes 
the rational vector space $A\otimes_{\mathbb Z}\mathbb Q$
and ${\rm cone}(w_i\, : i\in I)$ is the convex cone 
of $A_{\mathbb Q}$ generated by the vectors 
$w_i$ indexed by $I$. 
A reference for the following 
result is~\cite{ADHL}*{Exer. 2.11}.

\begin{proposition}
\label{pro:poly}
Let $\mathbb C[x_1,\dots,x_r]$
be an $A$-graded polynomial ring 
with homogeneous variables  
and let $w_i := \deg(x_i)$ for all $i$. This polynomial ring 
is a {\em Cox ring} 
if the following conditions hold:
\begin{enumerate}
\item
$A$ is generated by any $r-1$ elements
of $\{w_1,\dots,w_r\}$;  
\item
for each $1\leq i,j\leq r$ the interior of the cones 
${\rm cone}(w_k\, : k\neq i)$ and 
${\rm cone}(w_k\, : k\neq j)$ have non-empty 
intersection.
\end{enumerate}
\end{proposition}

The first condition is equivalent
to the triviality of the stabilizer of a general point
of the divisor $V(x_i)$ for all $i$. 
The second condition is used to guarantee 
that the quotient is good. As an application,   
the first algebra in Example~\ref{ex-1} is
a Cox ring, while the second one is not.

\section{Toric varieties}   \label{section.toric}

Proposition~\ref{pro:poly}
characterizes the graded polynomial rings 
 that are Cox rings. Now we would like to identify the varieties with such polynomial Cox rings.     
These are the {\em toric varieties} and we describe them in this section.    
The~affine variety $T := (\mathbb C^*)^n$ is called the $n$-dimensional algebraic torus or simply the torus.
The function $T \times T \rightarrow T$ given by coordinatewise multiplication is a morphism and it makes a torus $T$ into an algebraic group. 
\begin{definition}
A \emph{toric variety} is a normal variety that contains a torus $T$ as an open subset, and which has an action of $T$ via a morphism $T \times X \rightarrow X$ that restricted to $T \subseteq X$ is the $T \times T \rightarrow T$ coordinatewise multiplication action of $T$ on itself.  
\end{definition}
\begin{example}
Projective space $\mathbb{P}^n$ is a toric variety with the $T$-action given by 
$t \cdot x =  [x_0:t_1x_1:\dotsb:t_nx_n]$
for any $t=(t_1,\ldots,t_n) \in T$ and $x=[x_0:x_1:\dotsb:x_n] \in \mathbb{P}^n$. 
The coordinate chart $U_0 := \mathbb P^n
\setminus V(x_0)$ is an affine 
space invariant under the same action,
so that $\mathbb A^n$ is also a toric
variety.
\end{example}
In order to provide a description
for the Cox ring of a toric variety
we need to review the 
combinatorial language   
used to describe such varieties.
To begin with, each one-parameter
subgroup $\mathbb C^*\to T$ is
uniquely determined by a vector
of exponents in a space 
$N \simeq\mathbb Z^n$.
Dually each character $T\to\mathbb C^*$
is determined by a vector $u$ in
the dual $M = {\rm Hom}(N,\mathbb{Z}) \simeq\mathbb Z^n$ of $N$. The duality is 
expressed by a unimodular bilinear
pairing $M\times N\to\mathbb Z$,
$(u,v)\mapsto\langle u,v\rangle$.
This extends to a duality of the corresponding
rational vector spaces $M_{\mathbb Q}=M\otimes_{\mathbb Z}\mathbb Q$,
$N_{\mathbb Q}=N\otimes_{\mathbb Z}\mathbb Q$. 
If $X$ is an affine toric variety and
$x_0\in X$ is a point in the open torus 
orbit, the set of one-parameter 
subgroups $t\mapsto t^v$ such 
that $\lim_{t\to 0}t^v\cdot x_0$
exists in $X$ form the set of integer points
of a convex polyhedral (i.e., finitely
generated) cone
$\sigma\subseteq N_{\mathbb Q}$.
This cone $\sigma$ is strictly convex (i.e., $\{0\}$ is a face of $\sigma$), 
since $T$ is an invariant open subset of $X$.  
On the other hand, given such a cone,
its dual $\sigma^\vee := \{u\in 
M_{\mathbb Q}\, :\, \langle u,v\rangle
\geq 0\text{ for any $v\in \sigma$}\}$
defines the monoid algebra
$\mathbb C[\sigma^{\vee} \cap M]$
which turns out to be finitely generated.
The corresponding affine variety 
$X_\sigma$ is toric isomorphic 
to $X$, with torus action induced by 
the homomorphism
$\mathbb C[\sigma^{\vee} \cap M]\to
\mathbb C[M]\otimes
\mathbb C[\sigma^{\vee} \cap M]$
defined by $\chi^m\mapsto \chi^m\otimes\chi^m$.

More generally toric varieties are constructed
by gluing affine toric varieties which share
a common torus $T$.
Any toric variety $X$ has an open cover by affine toric varieties, and thus it corresponds to a 
collection $\Sigma$ of strictly convex 
polyhedral cones in the space $N_{\mathbb Q}$
and the gluing condition implies that
$\Sigma$ is a {\em fan}: 
each face of a cone in $\Sigma$ is again in 
$\Sigma$ and the intersection of any two cones in $\Sigma$ is a face of both.
If $\tau$ is a face of $\sigma$, 
and $m\in M$ is in the
relative interior of $\sigma^\vee\cap\tau^\perp$
then we have $\tau^\vee\cap M = (\sigma^\vee\cap M)
\oplus\mathbb Z_{\geq 0}\cdot (-m)$.
Equivalently the algebra
$\mathbb C[\tau^\vee\cap M]$ is the
localization of $\mathbb C[\sigma^\vee\cap M]$
at the multiplicative system
$S = \{1,\chi^m,\chi^{2m},\dots\}$,
which implies that $X_\tau$ is 
an open affine torus invariant subset 
of $X_\sigma$.
Using this, one can show that the affine toric varieties $U_{\sigma}$ for $\sigma$ in a fan $\Sigma$ can be glued to obtain a toric variety which we denote by $X_\Sigma$. 
If the fan $\Sigma$ arises from $X$ as above, then $X$ is isomorphic to $X_{\Sigma}$ as toric varieties, hence every toric variety arises from a fan.   

\begin{remark}
A toric variety $X$ has 
finitely many $T$-orbits and they 
correspond to the cones in the fan. 
For all the integer points 
$v$ in the relative interior of a cone $\sigma\in\Sigma$
the limits $\lim_{t\to 0}t^v\cdot x_0$ are in
the same torus orbit.
The toric variety $X$ has finitely many $T$-invariant subvarieties, 
each arising as the closure of a unique $T$-orbit. 
The $T$-invariant subvarieties are in inclusion-reversing correspondence with the cones in the fan. 
For example, $T$-invariant prime divisors
correspond to the rays of the fan (i.e., to the one-dimensional cones). 
We  denote the set of rays of $\Sigma$ 
by $\Sigma(1)$ and, by abuse of notation,
 use the same symbol for the set of
primitive generators in $N$ of the rays.
\end{remark}

Since the Laurent 
polynomial ring is factorial, the
divisor class group of a torus $T$ is
trivial
so that
any divisor of a toric
variety is linearly equivalent 
to a $T$-invariant one.
As a consequence, if $K$ is the 
group of $T$-invariant Weil divisors 
then the class map ${\rm cl}|_K\colon
K\to {\rm Cl}(X)$ is surjective.
Its kernel $K^0$ consists of $T$-invariant
principal divisors. Each such 
divisor $\div(\chi)$ cannot intersect the torus 
$T$, so that $\chi$ is a regular function on
$T$. Moreover $T$-invariance implies 
that $\chi|_T$ is a character of the torus. 
Vice versa, any character of $T$ gives a 
principal torus invariant divisor.
From now on we also assume that the only global regular invertible functions of $X$ are constants, which is equivalent to the rays of the fan $\Sigma$ spanning $N_{\mathbb{Q}}$. 
It follows that the map $M\to K^0$
defined by $u\mapsto\div(\chi^u)$
is an isomorphism of abelian groups.
Now, let $\mathcal S$ be as in~\eqref{defS}
and let $D_1,\dots,D_r$ be the prime 
$T$-invariant divisors of $X$.
The algebra $\mathbb C[x_1,\dots,x_r]
\otimes_{\mathbb C}\mathbb C[M]$
is $K$-graded by $\deg(x_i) := D_i$
and $\deg(\chi^m) := \div(\chi^{-m})$.
The homomorphism of 
$K$-graded algebras
\[
 \mathbb C[x_1,\dots,x_r]
 \otimes_{\mathbb C}\mathbb C[M]
 \to\Gamma(X,\mathcal S)
\]
defined by $x_i\otimes 1\mapsto 1 \in 
\Gamma(X,\mathcal S)_{D_i}$
and $1\otimes\chi^m\mapsto 
\chi^m\in\Gamma(X,\mathcal S)_{-\div(\chi^m)}$,
is an isomorphism.
The surjectivity follows from the
fact that for any $T$-invariant divisor $D$
the vector space $\Gamma(X,\mathcal S)_D$ 
is generated by monomials of $\mathbb C[M]$.
To prove the injectivity, first of all one shows
that $\div(\chi^m) = \sum_i\langle m,v_i\rangle D_i$,
where each $v_i$ is the primitive generator of the one-dimensional cone of the fan corresponding to $D_i$. 
Thus, if $D = \sum_ia_iD_i$, the degree $D$ 
part of each of the above algebras is isomorphic
to the vector space 
\[
 \bigoplus_{m \in \Delta_D\cap M} \mathbb C\cdot\chi^m,
\]
where $\Delta_D := \{m \in M_{\mathbb{Q}} 
\, :\, \langle m, v_i \rangle + a_i  \geq 0 \}$.
By the definition of Cox ring, $\mathcal R(X)$ 
is isomorphic to the quotient 
$\Gamma(X,\mathcal S)/\Gamma(X,\mathcal I)$,
where $\mathcal I$ is the ideal sheaf locally 
generated by $1-\chi^m$ for all $m\in M$.
As a consequence we have the following.
\begin{theorem}
The Cox ring of a toric variety $X$ 
with torus invariant prime divisors
$D_1,\dots,D_r$ is the
${\rm Cl}(X)$-graded polynomial
ring $\mathbb C[x_1,\dots,x_r]$, where 
$\deg(x_i) = [D_i]\in {\rm Cl}(X)$ 
for each $i$.
\end{theorem}

Next, we describe the correspondence
between full-dimensional lattice polytopes 
and pairs consisting of a projective toric variety together with an ample divisor.
A {\em lattice polytope}
$\Delta\subseteq M_{\mathbb Q}$
is a polytope with vertices in $M$.
Each maximal proper face $F$ of $\Delta$
corresponds to the generator $v_F\in N$
of the monoid of points of $N$ which 
lie on the inward normal to $F$.
Now, each face $Q$ of $\Delta$
defines a cone $\sigma_Q$ of $N_{\mathbb Q}$
generated by all the $v_F$ as $F$ runs
over the maximal proper faces of $\Delta$ 
which contain $Q$. The set of such cones
is the {\em normal fan} $\Sigma_\Delta$ to 
$\Delta$ and the corresponding toric 
variety is denoted by $X_\Delta$.
The ample divisor is 
$H 
 := 
 -\sum_{F}
 \min_{u\in\Delta}\langle u,v_F\rangle D_F,
$
where $F$ runs over all the maximal
faces of $\Delta$ and $D_F$ denotes 
the prime divisor of $X_\Delta$ defined
by the one-dimensional cone of the fan
generated by $v_F$.  
%Observe that $H$ is ample but not necessarily very ample, as can be seen when $\Delta$ is the three-dimensional simplex with vertices $(0,0,0),(1,0,0),(0,1,0),(1,1,3)$. In this case $H = 3D_F$, where $F$ is the facet generated by the last three vectors. The space of global sections of this divisor is generated by the characters $\chi^v$ where $v$ varies along the vertices of $\Delta$, so that the map that it defines is not injective on the torus $T$, being there of degree $3$.
Every pair consisting of a projective 
toric variety $X$ together with an ample 
$T$-invariant divisor $H$ of $X$ arises 
this way, where the polytope is $\Delta_H$
as defined above.
\begin{example}
Given at least two positive integers 
$a_0,\ldots,a_n$, such that any $n$
of them have greatest common divisor 
one, the {\em weighted projective space}
$\mathbb{P}(a_0,\ldots,a_{n})$
is the toric variety obtained as the good 
quotient of $\mathbb C^{n+1}\setminus\{0\}$
by the $\mathbb C^*$-action
$t\cdot (x_0,\dots,x_n) := 
(t^{a_0}x_0,\dots,t^{a_n}x_n)$.
Alternatively, if $l$ is the least
common multiple of the $n+1$
numbers, then $\mathbb{P}(a_0,\ldots,a_{n})$
is the toric variety defined by the lattice polytope
$\Delta\subseteq\mathbb Q^{n+1}$
with vertices $(l/a_0,0,\ldots,0)$, 
$(0,l/a_1,\ldots,0)$, \ldots, $(0,\ldots,0,l/a_n)$
contained in the hyperplane with equation
$a_0x_{0}+\cdots+a_{n}x_{n}=l$.
The divisor class group of the weighted 
projective space is free of rank one and
the Cox ring is the polynomial ring 
$\mathbb K[x_0, \ldots , x_n]$, 
graded by $\deg(x_i) := a_i$.
\end{example}

\subsubsection*{Cox rings of $T$-varieties}
We conclude the section by discussing 
the Cox ring of a $T$-variety. 

\begin{definition}
A \emph{$T$-variety} is a normal variety that 
admits an effective action of a torus $T$.
The difference $\dim(X) - \dim(T)$ is the 
{\em complexity} of $X$, so that the
$T$-varieties of complexity zero are the toric varieties.
\end{definition}

As before we assume that $X$ has only 
constant invertible global functions 
and finitely generated divisor class 
group.
For any point $x\in X$, let $T_x\subseteq T$ 
denote its stabilizer. 
The complement of the open 
subset $X_0 := \{x\in X\, :\, T_x\text{ is finite}\}$
is union of a finite number of prime divisors
$E_1,\dots,E_m$. According to a result of 
Sumihiro~\cite{sumihiro1974equivariant}*{Cor. 3}, there is a geometric 
quotient $q \colon X_0 \to X_0/T$ 
with an irreducible normal but possibly 
non-separated orbit space $X_0/T$. 
As shown in~\cite{HS}*{Proof of Thm. 1.2}
the prevariety $X_0/T$ admits a {\em separation}, that is, a 
surjective rational map 
$\pi\colon X_0/T\dashrightarrow Y$ 
onto a normal variety $Y$,
with only constant invertible global functions 
and a finitely generated divisor class group,
which is a local isomorphism over a big open
subset of $Y$. The same reference 
shows that there are prime divisors
$C_0,\dots,C_r$ of $Y$
such that $\pi$ is an isomorphism outside 
their union, that each $\pi^{-1}(C_i)$ 
is a disjoint union of prime divisors $C_{ij}$ 
and that all divisors with non-trivial finite
generic isotropy occur among the closures
$D_{ij}$ of the preimages of the $C_{ij}$ via
$q$. 
For each $i$ one defines a monomial
$x_i^{l_i} := \prod_jx_{ij}^{l_{ij}}$
where $l_{ij}$ is the order of the isotropy
group of a general point of $D_{ij}$.
If one denotes by $1_{C_i}$ a defining
section for $C_i$ in Cox coordinates on
$Y$, then, by~\cite{HS}*{Thm. 1.2}, 
the Cox ring of $X$ is isomorphic to
\[
\mathcal R(Y)[x_{ij},y_k]/
\langle x_i^{l_i} - 1_{C_i}\, :\, 0\leq i\leq r  \rangle,
\] 
where each $x_{ij}$ is a variable corresponding 
to a defining section of $D_{ij}$ and each 
$y_k$ is a variable corresponding to a defining section of $E_k$,
see~\cite{ADHL}*{Thm. 4.4.1.3}. 
When the complexity of $X$ is zero,  
we again see that the Cox ring is a 
polynomial ring. Indeed, in this case only
the variables $y_1,\dots,y_m$ survive
because there are no points with non-trivial
finite isotropy.
If $X$ has complexity one then $Y$ is the 
projective line $\mathbb P^1$ because of 
the finiteness of its divisor class group, 
and we see that in this case the Cox ring of $X$ is finitely generated. 
More generally, the Cox ring of $X$ is finitely generated if and only if the Cox ring of $Y$ is finitely generated. 
%In this case $C_0,\dots,C_r$ are points of $\mathbb P^1$. Choosing representatives $\tilde C_i\in \mathbb C^2$, any three such vectors satisfy a linear relation which, after replacing  $\tilde C_i$ with the monomial $x_i^{l_i}$, becomes a trinomial relation in $\mathbb C[x_{ij},y_k]$ in the first group of variables. 

\section{Finite generation}   
We now discuss the finite generation
of Cox rings from a geometric
perspective. 
Given an $A$-graded algebra $R$ 
and a submonoid $H\subseteq A$, one
can form the corresponding 
{\em Veronese subalgebra}
\[
 R_H := \bigoplus_{a\in H} R_a
\]
which inherits finite generation from
$R$, if $H$ is finitely generated. To see this, fix a graded surjection
$\pi\colon\mathbb C[x_1,\dots,x_r] \to R$,
where $x_i$ is homogeneous of degree
$a_i := \deg(\pi(x_i))$ and let 
$Q\colon \mathbb Z^r\to A$ be the
homomorphism defined by 
$e_i\mapsto a_i$. 
The monoid 
$M := Q^{-1}(H)\cap\mathbb Z^r_{\geq 0}$
is finitely generated, being the intersection of 
two finitely generated submonoids of $\mathbb{Z}^r$~\cite{ADHL}*{Prop. 1.1.2.2}.
Thus the monoid algebra $\mathbb C[M]$ is finitely generated as well and $\pi$ maps
it surjectively onto $R_H$.
More generally, given a homomorphism 
$\phi\colon \mathbb Z^s\to A$,
the $\mathbb Z^s$-graded algebra
\[
 S := \bigoplus_{u\in\mathbb Z^s}R_{\phi(u)}
\]
is finitely generated.  
Indeed, the Veronese subalgebra $R_H$, where $H$ is the image of $\phi$, is finitely generated.
Let $\pi: \mathbb C[u_1,\dots,u_k]\to R_H$
be a surjection of graded algebras, where
each $u_i$ is homogeneous of degree 
$b_i := \deg(\pi(u_i))$. 
Let $P\colon\mathbb Z^n\to\mathbb Z^s$ 
be the kernel of $\phi$ 
and let $\eta\colon H\to \mathbb Z^s$
be a map of sets such that 
$\phi\circ \eta = {\rm id}$.
The Laurent polynomial ring
$\mathbb C[u_1,\dots,u_k]
\otimes_{\mathbb C}
\mathbb C[\mathbb Z^n]$
is graded by $\mathbb Z^s$
by giving degree $\eta(b_i)$ 
to the variable $u_i$ and degree 
$P(m)$ to $\chi^m$.
The homomorphism 
$\mathbb C[u_1,\dots,u_k]
\otimes_{\mathbb C}
\mathbb C[\mathbb Z^n]\to S$,
defined by 
$u_i\otimes 1\mapsto \pi(u_i)\in S_{\eta(b_i)}$ and
$1\otimes\chi^m\mapsto 1\in S_{P(m)}$, is a surjection 
of $\mathbb Z^s$-graded algebras.
Observe that if one replaces $\phi$ 
with a homomorphism of monoids 
then the previous arguments still prove 
finite generation of the corresponding 
graded algebras.
An immediate consequence of these observations
is that the finite generation of the Cox ring 
is equivalent to that of the algebra
\[
 \bigoplus_{(m_1,\ldots,m_n) \in \mathbb{Z}^n }
 \Gamma(X,O_X(m_1D_1+\cdots+m_nD_n))
\]
whenever the classes of $D_1,\dots,D_n$ 
generate ${\rm Cl}(X)$.
In fact, the above equivalence still holds when the classes of $D_1,\dots,D_n$ generate 
a subgroup of finite index of 
${\rm Cl}(X)$, see~\cite{ADHL}*{Prop. 1.1.2.5}.
These algebras are frequently used when studying the finite generation of Cox rings because one can choose the $D_i$ 
according to the geometry of the variety. Some authors also call such an algebra a Cox ring of the variety. 
Here we follow the convention that the Cox ring of a variety is the algebra graded over the whole divisor class group as defined in Section~\ref{section.graded.algebras}. Our terminology agrees with the original presentation of the Cox ring of a toric variety as a polynomial ring in \cite{Cox}*{Section 1}.

The section ring $R(X,D_1,\dots,D_r)$ of a collection of Weil divisors $D_1,\dots,D_r$ on a variety $X$ is defined as 
\[
 \bigoplus_{(m_1,\ldots,m_r) \in (\mathbb{Z}_{\geq 0})^r } \Gamma(X,\mathcal O_X(m_1D_1+\cdots+m_rD_r)).
\]

\begin{example}    \label{example.toric.cone}
We now give a combinatorial argument for the finite
generation of the section ring 
$R(X,F_1,\dots,F_r)$ 
for any 
Weil divisors $F_1,\dots,F_r$ 
on a $\mathbb{Q}$-factorial toric variety.
Let us first consider the case $r=1$
and let $D := F_1$.
Let us see that its section ring 
$R(X,D)$ 
is finitely generated.  
We may assume that $D$ is $T$-invariant and write $D=\sum a_{i}D_{i}$ where $D_1,\ldots,D_r$ are the torus invariant prime divisors and $a_1,\ldots,a_r$ are integers.  
Since $\Delta_{mD}=m\Delta_D$ for any positive integer $m$, the desired finite generation is equivalent to the finite generation of the monoid of integral points in 
$M_{\mathbb Q}\oplus \mathbb{Q}$ which are inside the rational polyhedral cone generated by $\Delta_D \times \{ 1 \}$ and $\Delta_{0} \times \{ 0 \}$. 
The monoid of integral points in a rational polyhedral cone is always finitely generated by Gordan's Lemma, then the finite generation of $R(X,D)$ follows. The figure below helps visualize the idea behind this finite generation. 
\begin{figure}[ht]
\begin{center}
\definecolor{ttffqq}{rgb}{0.2,1,0}
\definecolor{qqzzff}{rgb}{0,0.6,1}
\begin{tikzpicture}[line cap=round,line join=round,>=triangle 45,x=1cm,y=1cm,scale=0.4]
\fill[line width=1pt,color=qqzzff,fill=qqzzff,fill opacity=0.04] (-1.897786427335731,2.032211282910059) -- (2.102213572664269,2.032211282910059) -- (2.8347899583669087,3.2531719257477953) -- (-1.1652100416330913,3.2531719257477953) -- cycle;
\fill[line width=1pt,color=ttffqq,fill=ttffqq,fill opacity=0.1] (0.3983,2.33406) -- (1.4361144390795235,2.5172039827039376) -- (1.008778214086317,2.9140161916262017) -- (0.06464937793123189,2.9779116211128334) -- (-0.5903433930364766,2.6660103016044006) -- cycle;
\fill[line width=1pt,color=qqzzff,fill=qqzzff,fill opacity=0.05] (-4.744466068339327,5.080528207275147) -- (5.255533931660672,5.080528207275147) -- (7.086974895917272,8.132929814369488) -- (-2.9130251040827284,8.132929814369488) -- cycle;
\fill[line width=1pt,color=ttffqq,fill=ttffqq,fill opacity=0.19] (0.99575,5.83515) -- (3.590286097698809,6.293009956759844) -- (2.5219455352157927,7.285040479065504) -- (0.16162344482807972,7.444779052782083) -- (-1.4758584825911916,6.665025754011001) -- cycle;
\draw [line width=1pt,color=qqzzff] (-1.897786427335731,2.032211282910059)-- (2.102213572664269,2.032211282910059);
\draw [line width=1pt,color=qqzzff] (2.102213572664269,2.032211282910059)-- (2.8347899583669087,3.2531719257477953);
\draw [line width=1pt,color=qqzzff] (2.8347899583669087,3.2531719257477953)-- (-1.1652100416330913,3.2531719257477953);
\draw [line width=1pt,color=qqzzff] (-1.1652100416330913,3.2531719257477953)-- (-1.897786427335731,2.032211282910059);
\draw [line width=1pt,color=ttffqq] (0.3983,2.33406)-- (1.4361144390795235,2.5172039827039376);
\draw [line width=1pt,color=ttffqq] (1.4361144390795235,2.5172039827039376)-- (1.008778214086317,2.9140161916262017);
\draw [line width=1pt,color=ttffqq] (1.008778214086317,2.9140161916262017)-- (0.06464937793123189,2.9779116211128334);
\draw [line width=1pt,color=ttffqq] (0.06464937793123189,2.9779116211128334)-- (-0.5903433930364766,2.6660103016044006);
\draw [line width=1pt,color=ttffqq] (-0.5903433930364766,2.6660103016044006)-- (0.3983,2.33406);
\draw [line width=1pt,color=qqzzff] (-4.744466068339327,5.080528207275147)-- (5.255533931660672,5.080528207275147);
\draw [line width=1pt,color=qqzzff] (5.255533931660672,5.080528207275147)-- (7.086974895917272,8.132929814369488);
\draw [line width=1pt,color=qqzzff] (7.086974895917272,8.132929814369488)-- (-2.9130251040827284,8.132929814369488);
\draw [line width=1pt,color=qqzzff] (-2.9130251040827284,8.132929814369488)-- (-4.744466068339327,5.080528207275147);
\draw [line width=1pt,color=ttffqq] (0.99575,5.83515)-- (3.590286097698809,6.293009956759844);
\draw [line width=1pt,color=ttffqq] (3.590286097698809,6.293009956759844)-- (2.5219455352157927,7.285040479065504);
\draw [line width=1pt,color=ttffqq] (2.5219455352157927,7.285040479065504)-- (0.16162344482807972,7.444779052782083);
\draw [line width=1pt,color=ttffqq] (0.16162344482807972,7.444779052782083)-- (-1.4758584825911916,6.665025754011001);
\draw [line width=1pt,color=ttffqq] (-1.4758584825911916,6.665025754011001)-- (0.99575,5.83515);
\draw [line width=0.5pt] (0,0)-- (-2.3613735721459066,10.664041206417602);
\draw [line width=0.5pt,densely dotted] (0,0)-- (0.25859751172492756,11.911646484451333);
\draw [line width=0.5pt,densely dotted] (0,0)-- (4.035112856345268,11.656064766504807);
\draw [line width=0.5pt] (0,0)-- (5.744457756318094,10.06881593081575);
\draw [line width=0.5pt] (0,0)-- (1.5932,9.33624);
\begin{scriptsize}
\draw [fill=black] (0,0) circle (1pt);
\draw [fill=black] (0.3983,2.33406) circle (1pt);
\draw [fill=black] (1.4361144390795235,2.5172039827039376) circle (1pt);
\draw[color=black] (3.2,2.5) node {$\Delta_{D} \times \{1\}$};
\draw [fill=black] (1.008778214086317,2.9140161916262017) circle (1pt);
\draw [fill=black] (0.06464937793123189,2.9779116211128334) circle (1pt);
\draw [fill=black] (-0.5903433930364766,2.6660103016044006) circle (1pt);
\draw [fill=black] (0.99575,5.83515) circle (1pt);
\draw [fill=black] (3.590286097698809,6.293009956759844) circle (1pt);
\draw[color=black] (5.8,6.3) node {$m\Delta_{D} \times \{m\}$};
\draw [fill=black] (2.5219455352157927,7.285040479065504) circle (1pt);
\draw [fill=black] (0.16162344482807972,7.444779052782083) circle (1pt);
\draw [fill=black] (-1.4758584825911916,6.665025754011001) circle (1pt);
\end{scriptsize}
\end{tikzpicture}
\caption{Visualizing the section ring of a divisor on a toric variety.}
\end{center}
\end{figure}
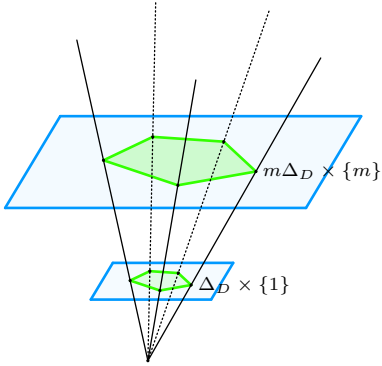

Now, for the general case, we may assume that each $F_j$ is a Cartier divisor and that it is a torus invariant divisor, as before. % of form $F_j=\sum a_{i,j}D_{i}$. 
In this case, the projectivization $\mathbb{P}(\mathcal{E})$ of the vector bundle $\mathcal{E}=\mathcal{O}_{X}(F_1) \oplus \cdots \oplus\mathcal{O}_{X}(F_r)$ is again a toric variety and the section ring  
$R(X,F_1,\dots,F_r)$ can be identified with the section ring  $R(\mathbb{P}(\mathcal{E}),H_{\mathcal E})$ where $-H_{\mathcal E}$ is a tautological divisor on $\mathbb{P}(\mathcal{E})$, and the restriction of $H_{\mathcal E}$ 
to each fiber of $\mathbb P(\mathcal E)$ is linearly equivalent to a hyperplane.
The finite generation of the section ring
$R(X,F_1,\dots,F_r)$
follows from the case $r=1$ considered before. 
\end{example}

\begin{example}[{Zariski's finite generation of section rings of semiample Cartier divisors}] 
If $D$ is a semiample Cartier divisor on a 
normal variety $X$, then $\psi_{mD}\colon X\to
\mathbb P^N$ is a morphism for some positive integer
$m$. Let $Y$ be the image of $\psi_{mD}$.
Then $mD$ is linearly equivalent to the pullback of a
hyperplane section $H\subseteq Y$ and moreover
$(\psi_{mD})_*\mathcal{O}_X = \mathcal{O}_Y$.
It follows that the algebra $R(X,mD)$ 
is isomorphic to $R(Y,H)$, which is finitely generated since $H$ is very ample on $Y$.
The latter property depends on the fact
that, given a hyperplane $H_0\subseteq \mathbb P^n$
with $H_0\cap Y = H$,
the algebra $R(\mathbb{P}^N,mH_0)$ is finitely generated and the sections 
$H^0(\mathbb{P}^N,mH_0)$ surject onto 
$H^0(Y,mH)$ for $m$ sufficiently large.
Then $R(X,D)$ is finitely generated as well since its Veronese subalgebra $R(X,mD)$ corresponding to a subgroup with finite index is finitely generated.
Zariski proves more generally that for any finite collection of semiample Cartier
divisors $D_1,\ldots, D_r$ on a normal projective variety $X$ their section ring
$R(X,D_1,\dots,D_r)$ is finitely generated.
We can see why this holds from what we have already proved.
Passing to a finite index Veronese subalgebra, we can assume that each of $D_1, \ldots, D_r$ is base point free. 
Let us denote by $\mathcal{E} := \mathcal{O}_{X}(D_1) \oplus \cdots \oplus\mathcal{O}_{X}(D_r)$ the direct
sum vector bundle and notice that 
$H_{\mathcal E}$ is base point free,
where $-H_{\mathcal E}$ is a
tautological divisor.
In particular $H_{\mathcal E}$ is semiample, and hence $R(\mathbb{P}(\mathcal{E}),H_{\mathcal E})$ is finitely generated. The claim now follows by noticing that the latter section ring is isomorphic to the section ring $R(X,D_1,\dots,D_r)$. 
\end{example}

\begin{example}
Fano varieties and more generally log-Fano varieties have finitely generated Cox rings, see \cite{birkar2010existence}.
\end{example}

\subsubsection*{Behavior of finite generation under morphisms}

How does the finite generation of Cox rings 
behave under maps?
This is a wide open question and 
here we just discuss it in some special
cases.
In what follows we  assume $X$ 
to be a normal variety that has only 
constant invertible global functions and a 
finitely generated divisor class group.
The algebras of sections over open subsets 
of $X$ of the sheaves of algebras considered 
before are finitely generated,  
if the Cox ring is so.  
More precisely in~\cite{Baker2011}*{Thm. 1.2}
it is proved that if the Cox ring of $X$ is finitely generated, then 
for any finitely generated 
subgroup $K$ of Weil divisors 
and any open subset $U \subseteq X$, the algebra of sections $R(U,K)$
is also finitely generated (see also~\cite{ADHL}*{Exer. 1.18}
for a description of the Cox
ring of an open subset in the case when the global invertible regular functions on $U$
are constants).
As a consequence, if $f\colon X\to Y$
is a birational contraction then the 
Cox ring of $Y$ is finitely generated
whenever that of $X$ is, being the former 
isomorphic to the Cox ring of the 
complement of the exceptional 
divisor.
More generally, even if $\mathcal R(X)$
is not finitely generated,
if $E$ is the exceptional
divisor (for simplicity assume it to be 
irreducible), then the pushforward
induces a surjection $\mathcal R(X)\to
\mathcal R(Y)$ with kernel generated by 
$x_E-1$, where $x_E$ is a defining section
for $E$ in Cox coordinates.
A more general result about finite generation
is the following.

\begin{theorem}[\cite{Ok}]
\label{theorem.Okawa}
If $f\colon X\to Y$ is a surjective morphism 
of normal projective $\mathbb Q$-factorial 
varieties and $X$ has finitely generated Cox 
ring, then the Cox ring of $Y$ is finitely generated.
\end{theorem}
When $f$ has connected fibers 
the statement is a consequence of the above 
observations on Veronese subalgebras.
Indeed, in this case $f_*\mathcal{O}_X=\mathcal{O}_Y$, and hence given a Cartier divisor $D$ on $Y$, by the projection formula 
there is an isomorphism
$\Gamma(X,\mathcal O_X(f^*D)) \simeq
\Gamma(Y,\mathcal O_Y(D))$.
It follows that if $K$ is a subgroup of
Cartier divisors of $Y$ whose image in
${\rm Cl}(Y)$ has finite index, 
there is an isomorphism between
$R(Y,K)$ and a Veronese subalgebra of the Cox ring of $X$. 
Thus $R(Y,K)$ is finitely generated 
and so the Cox ring of $Y$ is also finitely 
generated as well.

The finite generation of the Cox ring is preserved by passing to a good quotient
as shown in~\cite{Baker2011}*{Thm. 1.1}.
More precisely,
if $X$ is a normal variety with finitely generated Cox ring,
$G$ is a reductive affine algebraic group acting on $X$ 
and $U \subseteq X$ is an open invariant subset 
admitting a good quotient  $f\colon U \to U /\!\!/ G$,
then the Cox ring of $U/\!\!/G $ is finitely generated 
if this space has only constant invertible 
global functions.

\section{Mori chambers}  \label{section.mori.chambers}
Let us go back to the question:
``when do two varieties have isomorphic 
Cox rings?''.
By Remark~\ref{remark.small.modification}
this is the case when the two varieties
are isomorphic in codimension one.
Here we show that for normal projective $\mathbb{Q}$-factorial varieties
with finitely generated Cox ring the converse
is true: if two of them share the same Cox ring 
then they are isomorphic in codimension one.
Moreover in this case the number of such 
varieties is finite.
To begin with, the {\em base locus} 
of a divisor $D$ on a projective variety $X$ 
is the intersection of the supports of all the 
effective divisors linearly equivalent to $D$.
The {\em stable base locus} of $D$, denoted $B(D)$ is the intersection of the base loci of all the positive multiples $nD$ of $D$.
The divisor is {\em movable} if $B(D)$ has codimension at least two.
Assume now that ${\rm Cl}(X)$ is finitely 
generated. 
Hence, the rational vector space ${\rm Cl}_{\mathbb Q}(X) = {\rm Cl}(X) \otimes_{\mathbb Z}\mathbb Q$ is finite dimensional and has a well-defined Euclidean topology. 
The {\em effective cone} ${\rm Eff}(X)$
is the convex cone of ${\rm Cl}_{\mathbb Q}(X)$ generated by the classes of effective divisors
and the {\em moving cone} ${\rm Mov}(X)$
is the convex cone generated by classes of movable 
divisors. 
The closure $\overline{\rm Eff}(X)$ of ${\rm Eff}(X)$ in ${\rm Cl}_{\mathbb Q}(X)$ 
is called the
{\em pseudoeffective cone} of $X$. 
It is clear from the definition that
\[
 {\rm Mov}(X)
 \subseteq {\rm Eff}(X)
 \subseteq \overline{\rm Eff}(X).
\]

\begin{remark}  \label{remark.effective.cone}
If $U\subseteq X$ is an open subset with
complement of codimension at least $2$, then
the restriction map induces an isomorphism 
of divisor class groups under which the 
(pseudo)effective cones and the moving 
cones are identified. As a consequence 
the same identifications hold for a birational 
map which is an isomorphism in
codimension one.
\end{remark}

\begin{remark}
\label{remark.fg.cone}
For normal projective varieties the finite generation of 
${\rm Mov}(X)$ implies that of 
${\rm Eff}(X)$, see
\cite{ADHL}*{Lem. 4.3.3.3}, 
and both cones are uniquely determined 
by the degrees of a minimal generating
set of the Cox ring, see
\cite{ADHL}*{Prop. 3.3.2.3}.
Moreover the polyhedrality of pseudoeffective 
cones is preserved by surjections between normal projective varieties, 
see~\cite{CLTU2020}*{Lem. 2.2}.
\end{remark}

\begin{figure}[ht]   
 \begin{center}
 \begin{tikzpicture}[scale=0.6]
 \fill[fill=gray,fill opacity=0.6] (0,0) -- (2.6,2.6) -- (1,3) -- cycle;
 \draw [-] (0,0) -- (3,0);
 \draw [-] (0,0) -- (-1,3);
 \draw [-] (0,0) -- (1,3);
 \draw [-] (0,0) -- (3,1);
\draw [densely dotted,-,thick] (0,0) -- (2.3,2.3);
  \foreach \x/\y in {0/0,1/0,0.33/1,1.5/0.5,-0.5/1.5,2/1.5}
   {
     %\draw [-, thick] (0,0) -- (\x,\y);
     \fill (\x,\y) circle (0.7mm);
    }
 \draw (1,0) node[below] {\footnotesize $w_r$};
 \draw (1.5,0.5) node[below right] {\footnotesize $w_k$};
 \draw (0.5,1.3) node[left] {\footnotesize $w_j$};
 \draw (1,1.6) node[above] {\footnotesize $w$};
 \draw (-0.5,1.5) node[left] {\footnotesize $w_1$};
 \draw (1,1.6) circle (0.6mm);
 \end{tikzpicture}
 \end{center}
 \caption{The effective cone and its Mori chamber decomposition. 
 When the Cox ring is finitely generated every nef divisor is semiample.} \label{figure.semiample.nef}
\end{figure}
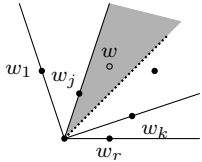

Assume now that $X$ is a normal projective $\mathbb{Q}$-factorial variety  
with a finitely generated Cox ring.
Given an effective Weil divisor $D$, the 
$\mathbb Z_{\geq 0}$-graded
algebra $R(X,D)$ is finitely generated so that it 
defines an affine variety equipped with a 
$\mathbb C^*$-action.
Fix a minimal set of homogeneous 
generators $g_1,\dots,g_s$ 
of degrees $d_1,\dots,d_s$ and define
the rational map
\[
 \varphi_D\colon X\dashrightarrow  
 \mathbb P(d_1,\dots,d_s)
 \quad
 p\mapsto [g_1(p):\cdots: g_s(p)], 
\] 
where the codomain is a weighted 
projective space and denote the 
closure of the image by $X(D)$.
\begin{definition}
Two effective divisors $D_1$, $D_2$ of $X$ 
are {\em Mori equivalent} if the following holds:
\begin{enumerate}
\item
$B(D_1) = B(D_2)$;
\item
there is an isomorphism 
$\phi\colon X(D_1)\to
X(D_2)$ such that the following diagram
commutes:
\[
\xymatrix{
 & X\ar@{-->}[ld]_-{\varphi_{D_1}}\ar@{-->}[rd]^-{\varphi_{D_2}}\\
 X(D_1)\ar[rr]^-\phi && X(D_2)
}
\]
\end{enumerate}
\end{definition}
Mori equivalence descends to an equivalence
relation on the effective cone. Each equivalence 
class is the relative interior of a rational polyhedral
cone called a {\em Mori chamber} of $X$.
The {\em nef cone} ${\rm Nef}(X)$, closure
of the cone of ample classes, is one of these
chambers. In particular the nef cone 
is generated by finitely many classes and
one can show that each nef class is semiample.
To give an idea of this fact assume for 
simplicity that $\Cl(X)$ has rank $2$. 
Let $w_1,\dots,w_r$ be the degrees
of a minimal set of generators $f_1,\dots,f_r$ of the Cox ring. 
Notice that $w_1,\dots,w_r$ span the effective cone ${\rm Eff}(X)=\overline{\rm Eff}(X)$, because every homogeneous
section $f\in\mathcal R(X)_w$ is a linear
combination of monomials
$\prod_if_i^{a_i}$
and thus $w = \sum_ia_iw_i$.
Moreover, by the combinatorial
description of Mori chambers (given 
below in this section) the nef cone of
$X$ is generated by two of these classes,
say $w_j$ and $w_k$, as in
Figure~\ref{figure.semiample.nef}. Within the nef cone there can be other
$w_i$ represented by black dots. 
We choose the darker region as in the figure,  
such that in its interior there are 
no degrees $w_i$ of the generators of the Cox ring. Each class
$w$ in this region is in the interior of the 
nef cone and thus is ample by Kleiman's
criterion.
Then, the stable base locus of $w$ 
is empty. Observe that $w$ is in the cone 
generated by $\{w_i\, :\, i\in I\}$ if and 
only if there are an integer $n>0$ and 
a monomial of degree $nw$ of the form
$\prod_{i\in I}f_i^{a_i}$. Since the stable
base locus $B(w)$ is the intersection of
the zero loci of these monomials we 
deduce that $B(w_j)\subseteq B(w)$
because each cone in the $w_i$ which
contains $w$ must contain also $w_j$.
Thus $B(w_j) = \emptyset$ or equivalently 
$w_j$ is semiample.

\begin{definition}
\label{mds}
A normal projective $\mathbb Q$-factorial
variety $X$ with a finitely generated divisor class group 
is a {\em Mori dream space}
if there are normal projective $\mathbb Q$-factorial varieties 
$X_1,\dots,X_n$ and isomorphisms in codimension one
$\phi_i\colon X\dashrightarrow X_i$, 
such that each nef cone $\Nef(X_i)$ is generated 
by finitely many semiample classes and 
\[ 
\Mov(X)
 = \bigcup_{i=1}^n\phi_i^*\Nef(X_i). 
\]
\end{definition}

The following theorem, one of the main
results of~\cite{HK}, shows that the
definition above fits perfectly into the context of
this note.

\begin{theorem}
Let $X$ be a normal projective 
$\mathbb Q$-factorial variety with finitely 
generated divisor class group. Then
$X$ is a Mori dream space if and only
if the Cox ring of $X$ is finitely generated.
\end{theorem}

As a consequence of the theorem, it makes 
sense to talk about Mori chambers of a 
Mori dream space. In this case one shows
that the cones in the decomposition of 
$\Mov(X)$ given in Definition~\ref{mds}
are all full-dimensional Mori chambers.
A combinatorial description of these chambers
in terms of a presentation of the Cox ring can
be given as explained below.
Given a subset $I\subseteq \{1,\dots,r\}$, 
denote by $C_I$ the cone of ${\rm Cl}_{\mathbb Q}(X)$
generated by $\{w_i\, :\, i\in I\}$.
The cone $C_I$ corresponds to the reduced
monomial $\prod_{i\in I}f_i$ in the generators
of the Cox ring.
It is not difficult to show that the moving cone
of $X$ is~\cite{ADHL}*{Prop. 3.3.2.9}
\[
 \Mov(X) = \bigcap_{i=1}^rC_{\{1,\dots,r\}\setminus\{i\}}.
\]
Indeed, for any class $w$ in the right-hand side
and any $i\in\{1,\dots,r\}$ there is an $n>0$ and 
a monomial in $\mathcal R(X)_{nw}$ which does 
not contain $f_i$. It is also possible to give a 
combinatorial description of Mori chambers:
the Mori chamber which contains the class
$w$ is
\[
 \lambda(w) = \bigcap_{w\in C_I}C_I.
\]
Finally observe that 
Mori equivalence refines the equivalence 
relation $D_1\sim D_2$ if $B(D_1) = B(D_2)$.
The equivalence classes for this relation
are the relative interiors of unions of 
some Mori chambers.
The two equivalence relations can be
distinct for a given Mori dream 
space, as shown by the following picture which displays the Mori chamber
decomposition of the effective cone 
of a toric threefold. The degrees of
generators of the Cox ring are the 
six black dots.

\begin{figure}[ht]
\begin{center}
\begin{tikzpicture}[scale=.3]
\tkzDefPoint(-3,0){P1}
\tkzDefPoint(3,0){P2}
\tkzDefPoint(6,6){P3}
\tkzDefPoint(-6,6){P4}
\tkzDefPoint(0,5){P5}
\tkzDefPoint(0,3.5){P6}
\tkzInterLL(P1,P5)(P2,P4) \tkzGetPoint{Q1}
\tkzInterLL(P2,P5)(P1,P3) \tkzGetPoint{Q2}
\tkzInterLL(P2,P4)(P1,P3) \tkzGetPoint{Q3}
\tkzInterLL(P1,P6)(P2,P4) \tkzGetPoint{Q4}
\tkzInterLL(P2,P6)(P1,P3) \tkzGetPoint{Q5}
\tkzInterLL(P4,P6)(P1,P5) \tkzGetPoint{Q6}
\tkzInterLL(P3,P6)(P2,P5) \tkzGetPoint{Q7}
\tkzFillPolygon[color = black](P6,Q4,Q3,Q5)
\tkzFillPolygon[color = gray!50](P5,P6,Q6)
\tkzFillPolygon[color = gray!50](P5,P6,Q7)
\tkzDrawSegments[thick](P1,P2 P2,P3 P3,P4 P4,P1)
\tkzDrawSegments[densely dotted](P5,P1 P5,P2 P5,P3 P5,P4 P6,P1 P6,P2 P6,P3 P6,P4 P5,P6 P4,P2 P3,P1)
\tkzDrawPoints[fill=black,color=black,size=3](P1,P2,P3,P4,P5,P6)
\end{tikzpicture}
\end{center}
 \caption{Mori chamber
decomposition of the effective cone 
of a toric threefold.}
\end{figure}
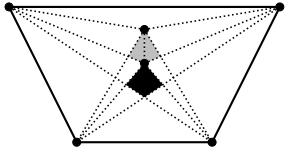

There are $17$ full-dimensional Mori 
chambers, five of 
which are inside the moving cone. The nef 
cone is the black chamber, while the two
gray chambers have the same stable 
base locus~\cite{LMR}.

\begin{remark}
We end-up with a remark about
GIT chambers for the action of the
quasitorus $H_X$ on the affine variety
$\overline X$. Given an effective class 
$w\in\mathbb {\rm Cl}_{\mathbb Q}(X)$
one defines the open subset
$\overline X^{\rm ss}(w)\subseteq\overline X$
consisting of the $\bar x\in \overline X$ 
such that $f(\bar x)\neq 0$ for some $n>0$
and $f\in\mathcal R(X)_{nw}$.
This is the subset of {\em semistable
points} defined by $w$, which is an
open $H_X$-invariant subset of
$\overline X$. The set of effective 
classes that give rise to the same
subset of semistable points is the 
relative interior of a polyhedral cone. 
These cones are called 
GIT chambers and one of the main
results of~\cite{HK} states that GIT
chambers are exactly Mori chambers.
\end{remark}

\section{Examples}

We  now illustrate some approaches that have been used in the literature to determine the finite or non-finite generation of the Cox rings of notable varieties. 
While the complete proofs frequently involve arguments specific to the geometry of the variety under consideration, we mention a few ideas with great applicability. We include references for the use of such ideas in concrete examples. 

Given a normal projective $\mathbb{Q}$-factorial variety with a finitely generated Cox ring (i.e., a Mori dream space), the cones of nef, pseudoeffective and movable divisors are rational polyhedral. Hence, one can disprove the finite generation of the Cox ring by studying these cones. 
For example, the pseudoeffective cone of the moduli space $\overline{M}_{1,n}$ 
of stable genus-one curves with $n$ ordered marked points 
is not polyhedral for $n \geq 3$.  
Hence, the moduli space $\overline{M}_{1,n}$ does not have a finitely generated Cox ring for $n \geq 3$ \cite{chen2014extremal}. 

\begin{example}
The finite or non-finite generation of the Cox ring has been decided for the
blow up of $\mathbb{P}^{r-1}$ at $n > r \geq 2$ 
points in very general position: 
finite generation occurs if and only if 
(see~\cite{mukai2004} and~\cite{castravet2006hilbert}*{Thm. 1.3})
\[
\frac{1}{r}+\frac{1}{n-r} > \frac{1}{2}. 
\]
If instead the set $S$ of points lies on a 
hyperplane $H\subseteq\mathbb P^n$, 
with $n>2$,  then, the Cox ring 
of $\operatorname{Bl}_S {\mathbb{P}}^n$ 
is isomorphic to a polynomial ring in one 
variable over the Cox ring of 
$\operatorname{Bl}_S H$,
see~\cite{GHPS}*{Lem. 4.3}.
%This last result can be proved using that $\operatorname{Bl}_S {\mathbb{P}}^n$ is a $T$-variety for $T = \mathbb C^*$.

\end{example}

Since nef divisors on a Mori dream space are semiample,   
one may establish that a Cox ring is non-finitely generated by exhibiting a nef divisor on the variety that is not semiample. 
This approach has been helpful when studying the Cox rings of blowups at a general point of weighted projective planes. 
These Cox rings are also studied in commutative algebra since their finite generation is equivalent to the Noetherianity of the extended saturated Rees algebra of a monomial prime ideal.
For a concrete example using this approach see Example~\ref{example.12.13.17}, which adapts a general argument from \cite{gonzalez2016some}.

\begin{example}
[$\operatorname{Bl}_{e}\mathbb{P}(12,13,17)$ is not a Mori dream space \cite{gonzalez2016some}]  
\label{example.12.13.17}
The triangle $\Delta$ with integral 
vertices $(11,-26)$, $(50,0)$, $(-1,34)$
defines a toric pair consisting of the 
weighted projective plane 
$\mathbb P_\Delta = \mathbb{P}(12,13,17)$
together with the ample class $H$
of self-intersection $H^2 = 52\cdot 51$
equal to twice the area of $\Delta$.
Let $\pi\colon X\to \mathbb P_{\Delta}$ be the blow up 
at the unit element of the torus $e\in\mathbb P_\Delta$ with 
exceptional divisor $E$.
The binomial $1-y$ defines a curve
in $\mathbb P_{\Delta}$ which passes 
simply through $e$ so that its strict transform 
$C\subseteq X$ is an irreducible
curve with $C\cdot E = 1$.
Notice that the triangle $\Delta$ supports the Laurent polynomial 
$f(x,y) := x^{11}y^{-26}(1-y)^{52}$, which has a nonzero coefficient
over each edge of $\Delta$, and which vanishes to order $52$ at $e$.  
Therefore the class of $C$
is $\frac{1}{52}\pi^*H-E$.
Using the above intersection numbers,
the self-intersection of $C$ on $X$ can be computed as
\[
 C^2=\left(\frac{1}{52^2}(\pi^*H)^2+E^2\right)=\frac{52\cdot 51}{52^2}-1<0.
\] 
Therefore, $C$ is an irreducible negative curve on $X$.
Let us show that $X$ is not a Mori dream space.
Since $\operatorname{Cl}(X)$ has rank two,  
the pseudoeffective cone is generated by the 
classes of the two negative curves $C$ and $E$.
Taking dual with respect to the intersection
form, we see that the nef cone of $X$ is generated by 
$\pi^*H$ together with any class $D$ satisfying $D \cdot E >0$ and $D \cdot C = 0$. 
We can then take $D := \pi^*H-51E$. 
Notice that $\Delta$ contains $(49,0)$ and $(50,0)$, but the remaining lattice points in $\Delta$ have an $x$-coordinate with value at most $48$.
Then, the partial derivative $\partial_x^{49} \partial_y |_{(x=1,y=1)}$ vanishes when evaluated on the monomial $x^a y^b$ for each lattice point in $\Delta$, except its left vertex $(-1,34)$.
It follows that any Laurent polynomial $g(x,y)$ supported on $\Delta$ that vanishes to order $51$ at $e$ has the coefficient corresponding to the left vertex of $\Delta$ equal to zero. 
Indeed, $\partial_x^{49} \partial_y g(x,y) |_{(x=1,y=1)}$ is equal to a nonzero multiple of the coefficient of $x^{-1}y^{34}$ in $g(x,y)$, but on the other hand it is equal to zero since $g(x,y)$ has order $51$ at $e$. 
Then every effective divisor linearly equivalent to $D$ passes through the point in $X$ that is mapped by $\pi$ to the torus invariant point of $X_\Delta$ corresponding to the left vertex $(-1,34)$ of $\Delta$. 
In particular, $D$ is not base point free. 
Given any positive integer $m$, we can translate the triangle $m\Delta$ such that the left vertex has the form $(-1,a)$ and the right vertex is $(51m-1,0)$, and we can repeat the same argument using now the partial derivative $\partial_x^{51m-2} \partial_y |_{(x=1,y=1)}$ to conclude that $mD$ is not base point free. 
Then, $D$ is a nef divisor that is not semiample on $X$, and therefore $X=\operatorname{Bl}_{e}\mathbb{P}(12,13,17)$ is not a Mori dream space. 
\end{example}

One of the most fruitful methods to decide the finite or non-finite generation of the Cox ring of one variety is to appropriately change the variety under consideration to another more suitable for the finite generation question. This is frequently done using the behavior of finite generation under isomorphisms in codimension one (Rmk.~\ref{remark.small.modification}
and Rmk.\ref{remark.effective.cone}) and surjective morphisms  (Thm.~\ref{theorem.Okawa}
and Rmk.~\ref{remark.fg.cone}), as follows. 
Suppose that there is a chain of rational maps 
\begin{equation}   \label{equation.chain}  
    X_{1} \dasharrow X_{2} \dasharrow \cdots \cdots \dasharrow X_{n-1} \dasharrow X_{n}
\end{equation}
where each $X_i$ is a normal projective $\mathbb{Q}$-factorial variety and each map is either a surjective morphism or an isomorphism in codimension one.
If $X_1$ has a finitely generated Cox ring then the same 
holds for $X_{n}$.
To illustrate these ideas we  now outline the construction of one such chain of rational maps, starting from the moduli space $\overline{M}_{0,n}$ for any $n \geq 10$ and ending on the blown-up toric surface $\operatorname{Bl}_{e} \mathbb{P}(12,13,17)$ discussed in Example~\ref{example.12.13.17}. This construction shows that the Cox ring of $\overline{M}_{0,n}$ is not finitely generated for $n \geq 10$.  

\subsubsection*{The moduli space $\overline{M}_{0,n}$}

The moduli space $M_{0,n}$ parameterizes configurations of $n$ distinct ordered points in $\mathbb{P}^1$ up to automorphisms. 
A family of such configurations may have as a limit a configuration where some of the points collide, and therefore $M_{0,n}$ is not compact.  
The group $\operatorname{Aut}(\mathbb{P}^1)=\operatorname{PGL(2)}$ acts on triples of distinct ordered points in $\mathbb{P}^1$, transitively and with trivial stabilizers. Hence, for $n \geq 4$ we have the identification 
\[
M_{0,n}=(\mathbb{P}^1\smallsetminus \{0,1,\infty\} )^{n-3} \smallsetminus \bigcup_{i \neq j} \Delta_{i,j},\] 
where $\Delta_{i,j}$ denotes the locus where the $i$-th and $j$-th components coincide.   
There is a well-known compactification of $M_{0,n}$ denoted by $\overline{M}_{0,n}$ and called the Grothendieck-Knudsen compactification of $M_{0,n}$.  
This compactification parameterizes proper connected curves with at worse nodal singularities, which are a union of 
irreducible components isomorphic to $\mathbb{P}^1$ with dual graph forming a tree, with $n$ ordered marked points in the smooth locus, and such that each component contains at least $3$ points that are either a node or a marked point. 
Kapranov showed that $\overline{M}_{0,n}$ is isomorphic to the iterated blow up in a suitable order of $\mathbb{P}^{n-3}$ along the strict transforms of all linear subspaces spanned by subsets of $n-1$ points in linearly general position. 

If we denote the $n-1$ points in this construction by $p_1,\ldots,p_{n-1} \in \mathbb{P}^{n-3}$ and we first perform the iterated blow up of $\mathbb{P}^{n-3}$ along the strict transforms of all linear subspaces spanned by $p_1,\ldots,p_{n-2}$ in order of increasing dimension, we obtain a toric variety called the Losev-Manin moduli space which we  denote by $LM_n$.  
We can assume that $p_{n-1}$ is the unit element $e$ of the torus of $LM_n$. 
If we now blow up $LM_n$ along the strict transforms of the remaining linear subspaces, all of them containing $p_{n-1}$, the resulting variety is  $\overline{M}_{0,n}$.
In particular, there exist surjective morphisms 
\[
\overline{M}_{0,n} \rightarrow \operatorname{Bl}_{e} LM_n \rightarrow LM_n.
\]
From this construction we also see that $\overline{M}_{0,n}$ is a nonsingular projective variety. 
For $n \leq 6$ the variety $M_{0,n}$ is log-Fano and hence it is a Mori dream space. 
Hu and Keel asked in \cite{HK} whether $\overline{M}_{0,n}$ was a Mori dream space in general. 
A negative answer for $n \geq 134$ was given in \cite{castravet2015overline}, with the bound improved to $n \geq 13$ in \cite{gonzalez2016some}, and later improved to $n \geq 10$ in \cite{hausen2018blowing}. 
Let us describe the main points of the proof that $\overline{M}_{0,n}$ is not a Mori dream space for $n \geq 10$. 
For each $n$ there exists a surjective morphism
\[
 \overline{M}_{0,n+1} \rightarrow \overline{M}_{0,n}
\]
called the forgetful morphism which intuitively forgets one of the marked points (for example the last one). 
Hence, given positive integers $n \geq m\geq 4$ there exists a surjective morphism $\overline{M}_{0,n} \rightarrow \operatorname{Bl}_e LM_m$. 
By Okawa's result on images of Mori dream spaces, it follows that if 
$\operatorname{Bl}_e LM_m$ is not a Mori dream space, 
then $\overline{M}_{0,n}$ is not a Mori dream space for all $n \geq m$.

\begin{example}  \label{CT.toric.projection}
Suppose that we have a morphism of projective $\mathbb{Q}$-factorial toric varieties $\mathbb P_{\Sigma'} \rightarrow \mathbb P_\Sigma$ induced by a surjective map of lattices $\pi\colon N' \rightarrow N$.
Suppose also that $\operatorname{ker}(\pi)$ is one-dimensional,
generated by two rays of $\Sigma'$, and that the rays of $\Sigma$ are precisely the images of the rays in $\Sigma'$.
Then the rational map $\operatorname{Bl}_{e} \mathbb P_{\Sigma'} \dashrightarrow \operatorname{Bl}_{e} \mathbb P_\Sigma$ factors as an isomorphism in codimension one followed by a surjection by~\cite{castravet2015overline}*{Prop. 3.1}.
Then, the codomain $\operatorname{Bl}_{e} \mathbb P_{\Sigma'}$ is a Mori
dream space if the domain $\operatorname{Bl}_{e} \mathbb P_{\Sigma}$ is so.
\end{example}

The rays in the fan of $LM_n$ are generated by the vectors $v$ in $N=\mathbb{Z}^{n-3}$ such that either $v$ or $-v$ has all entries in $\{0,1\}$. 
Let $N' \subseteq N$ be a saturated sublattice of rank $n-5$, generated by a subset of $\{0, 1\}^{n-3}$. Assume that there exist $v_1, v_2, v_3 \in \{0, 1\}^{n-3}$ whose images under the quotient map $\pi: N \rightarrow N/N'$ generate $N/N'\cong \mathbb{Z}^2$, and such that there are pairwise coprime positive integers $a, b, c$ such that $av_1 \pm bv_2 \pm cv_3 \in N'$. Applying Example~\ref{CT.toric.projection} iteratively, it follows that if $\operatorname{Bl}_e LM_n$ is a Mori dream space, then $\operatorname{Bl}_{e}\mathbb{P}(a,b,c)$ is also a Mori dream space. 
We can apply this to the homomorphism 
$\pi\colon N=\mathbb{Z}^7 \rightarrow \mathbb{Z}^2$ given by the matrix 
\[
\begin{bmatrix*}[r]
1 & 0 & 1 & -2 & -1 & 1  & 0 \\
0 & 1 & -1 & -3 & -2 & 2  & 1
\end{bmatrix*},
\]
with $N'=\operatorname{ker}(\pi)$
and with $v_1=e_1 + e_3 + e_4 + e_5$, $v_2=-(e_4 + e_5)$ and $v_3=-(e_1 + e_3 + e_6)$, and $a=12, b=13, c=17$. 
This shows that for each $n \geq 10$ there exists a chain of maps 
\[
\overline{M}_{0,n}=X_{0} \dasharrow  X_{1} \dasharrow \cdots  \dasharrow X_{r}= \operatorname{Bl}_{e} \mathbb{P}(12,13,17)
\]
such that all varieties are normal, projective and $\mathbb{Q}$-factorial, and such that each map is either a surjective morphism or an isomorphism in codimension one.
In Example~\ref{example.12.13.17} we showed that $\operatorname{Bl}_{e} \mathbb{P}(12,13,17)$ is not a Mori dream space, and therefore 
$\overline{M}_{0,n}$ is not a Mori dream space for each $n \geq 10$.

One can use the same projection 
$\pi\colon\mathbb Z^7\to\mathbb Z^2$
to prove that $\overline{\rm Eff}(\overline{M}_{0,n})$ 
is not polyhedral for all $n \geq 10$, as done 
in~\cite{CLTU2020}*{Thm. 1.3}. 
Indeed, one gets a chain of maps as in
 (\ref{equation.chain}), 
\[
\overline{M}_{0,n}=X_{0} \dasharrow  X_{1} \cdots  \dasharrow X_{r}= \operatorname{Bl}_{e} \mathbb{P}_\Delta, 
\]
where the images of the primitive generators of the rays of $LM_{10}$ via $\pi$
are points (both black and white) in the following
figure.
\begin{figure}[ht]   
\begin{center}
\begin{tikzpicture}[scale=.35, rotate= 90]  
\tkzInit[xmin=-3.5,xmax=3.5,ymin=-6.5,ymax=6.5]\tkzGrid
  \draw[line width=.2mm ] (-4,0) -- (4,0);
  \draw[line width=.2mm]  (0,-7) -- (0,7);

 \tkzDefPoint(-3,-5){P1}
 \tkzDefPoint(-3,-1){P2}
 \tkzDefPoint(-2,-6){P3}
 \tkzDefPoint(-1,-6){P4}
 \tkzDefPoint(-1,3){P5}
 \tkzDefPoint(1,-3){P6}
 \tkzDefPoint(1,6){P7}
 \tkzDefPoint(2,6){P8}
 \tkzDefPoint(3,1){P9}
 \tkzDefPoint(3,5){P10}

 \tkzDefPoint(1,-1){Q1}
 \tkzDrawPoints[fill=black,color=black,size=3](Q1)
 \tkzDefPoint(-1,0){Q2}
 \tkzDrawPoints[fill=black,color=black,size=3](Q2)
 \tkzDefPoint(1,-2){Q3}
 \tkzDrawPoints[fill=black,color=black,size=3](Q3)
 \tkzDefPoint(0,0){Q4}
 \tkzDrawPoints[fill=black,color=black,size=3](Q4)
 \tkzDefPoint(-1,1){Q5}
 \tkzDrawPoints[fill=black,color=black,size=3](Q5)
 \tkzDefPoint(1,-3){Q6}
 \tkzDrawPoints[fill=black,color=black,size=3](Q6)
 \tkzDefPoint(-1,2){Q7}
 \tkzDrawPoints[fill=white,color=white,size=3](Q7)
 \tkzDefPoint(0,1){Q8}
 \tkzDrawPoints[fill=black,color=black,size=3](Q8)
 \tkzDefPoint(-1,3){Q9}
 \tkzDrawPoints[fill=white,color=white,size=3](Q9)
 \tkzDefPoint(0,2){Q10}
 \tkzDrawPoints[fill=black,color=black,size=3](Q10)
 \tkzDefPoint(0,3){Q11}
 \tkzDrawPoints[fill=black,color=black,size=3](Q11)
 \tkzDefPoint(0,4){Q12}
 \tkzDrawPoints[fill=black,color=black,size=3](Q12)
 \tkzDefPoint(0,-1){Q13}
 \tkzDrawPoints[fill=black,color=black,size=3](Q13)
 \tkzDefPoint(-2,0){Q14}
 \tkzDrawPoints[fill=black,color=black,size=3](Q14)
 \tkzDefPoint(0,-2){Q15}
 \tkzDrawPoints[fill=black,color=black,size=3](Q15)
 \tkzDefPoint(-2,1){Q16}
 \tkzDrawPoints[fill=black,color=black,size=3](Q16)
 \tkzDefPoint(1,0){Q17}
 \tkzDrawPoints[fill=white,color=white,size=3](Q17)
 \tkzDefPoint(0,-3){Q18}
 \tkzDrawPoints[fill=black,color=black,size=3](Q18)
 \tkzDefPoint(1,1){Q19}
 \tkzDrawPoints[fill=black,color=black,size=3](Q19)
 \tkzDefPoint(0,-4){Q20}
 \tkzDrawPoints[fill=black,color=black,size=3](Q20)
 \tkzDefPoint(1,2){Q21}
 \tkzDrawPoints[fill=black,color=black,size=3](Q21)
 \tkzDefPoint(1,3){Q22}
 \tkzDrawPoints[fill=black,color=black,size=3](Q22)
 \tkzDefPoint(1,4){Q23}
 \tkzDrawPoints[fill=black,color=black,size=3](Q23)
 \tkzDefPoint(1,5){Q24}
 \tkzDrawPoints[fill=black,color=black,size=3](Q24)
 \tkzDefPoint(1,6){Q25}
 \tkzDrawPoints[fill=black,color=black,size=3](Q25)
 \tkzDefPoint(2,0){Q26}
 \tkzDrawPoints[fill=black,color=black,size=3](Q26)
 \tkzDefPoint(2,1){Q27}
 \tkzDrawPoints[fill=black,color=black,size=3](Q27)
 \tkzDefPoint(2,2){Q28}
 \tkzDrawPoints[fill=black,color=black,size=3](Q28)
 \tkzDefPoint(2,3){Q29}
 \tkzDrawPoints[fill=black,color=black,size=3](Q29)
 \tkzDefPoint(2,4){Q30}
 \tkzDrawPoints[fill=black,color=black,size=3](Q30)
 \tkzDefPoint(2,5){Q31}
 \tkzDrawPoints[fill=black,color=black,size=3](Q31)
 \tkzDefPoint(2,6){Q32}
 \tkzDrawPoints[fill=black,color=black,size=3](Q32)
 \tkzDefPoint(3,1){Q33}
 \tkzDrawPoints[fill=white,color=white,size=3](Q33)
 \tkzDefPoint(-3,-1){Q34}
 \tkzDrawPoints[fill=black,color=black,size=3](Q34)
 \tkzDefPoint(3,2){Q35}
 \tkzDrawPoints[fill=white,color=white,size=3](Q35)
 \tkzDefPoint(-3,-2){Q36}
 \tkzDrawPoints[fill=black,color=black,size=3](Q36)
 \tkzDefPoint(3,3){Q37}
 \tkzDrawPoints[fill=black,color=black,size=3](Q37)
 \tkzDefPoint(-3,-3){Q38}
 \tkzDrawPoints[fill=black,color=black,size=3](Q38)
 \tkzDefPoint(3,4){Q39}
 \tkzDrawPoints[fill=black,color=black,size=3](Q39)
 \tkzDefPoint(-3,-4){Q40}
 \tkzDrawPoints[fill=black,color=black,size=3](Q40)
 \tkzDefPoint(3,5){Q41}
 \tkzDrawPoints[fill=black,color=black,size=3](Q41)
 \tkzDefPoint(-3,-5){Q42}
 \tkzDrawPoints[fill=white,color=white,size=3](Q42)
 \tkzDefPoint(-2,-1){Q43}
 \tkzDrawPoints[fill=black,color=black,size=3](Q43)
 \tkzDefPoint(-2,-2){Q44}
 \tkzDrawPoints[fill=black,color=black,size=3](Q44)
 \tkzDefPoint(-2,-3){Q45}
 \tkzDrawPoints[fill=white,color=white,size=3](Q45)
 \tkzDefPoint(-2,-4){Q46}
 \tkzDrawPoints[fill=black,color=black,size=3](Q46)
 \tkzDefPoint(-2,-5){Q47}
 \tkzDrawPoints[fill=black,color=black,size=3](Q47)
 \tkzDefPoint(-2,-6){Q48}
 \tkzDrawPoints[fill=black,color=black,size=3](Q48)
 \tkzDefPoint(-1,-1){Q49}
 \tkzDrawPoints[fill=black,color=black,size=3](Q49)
 \tkzDefPoint(-1,-2){Q50}
 \tkzDrawPoints[fill=black,color=black,size=3](Q50)
 \tkzDefPoint(-1,-3){Q51}
 \tkzDrawPoints[fill=black,color=black,size=3](Q51)
 \tkzDefPoint(-1,-4){Q52}
 \tkzDrawPoints[fill=black,color=black,size=3](Q52)
 \tkzDefPoint(-1,-5){Q53}
 \tkzDrawPoints[fill=black,color=black,size=3](Q53)
 \tkzDefPoint(-1,-6){Q54}
 \tkzDrawPoints[fill=black,color=black,size=3](Q54)
 \tkzDefPoint(2,-1){Q55}
 \tkzDrawPoints[fill=black,color=black,size=3](Q55)

 \draw (-1,2) circle (1.5mm); %Q7
 \draw (-1,3) circle (1.5mm); %Q9
 \draw (1,0) circle (1.5mm); %Q17
 \draw (3,1) circle (1.5mm); %Q33
 \draw (3,2) circle (1.5mm); %Q35
 \draw (-3,-5) circle (1.5mm); %Q42
 \draw (-2,-3) circle (1.5mm); %Q45

\end{tikzpicture}
\end{center}
\caption{Images of the primitive generators of the rays of $LM_{10}$ under the projection $\pi$.}   \label{figure.points}
\end{figure}
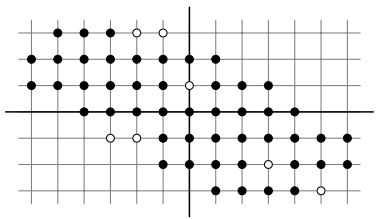

The white points in Figure~\ref{figure.points} generate the rays  
of the normal fan of the polygon
$\Delta$ whose vertices are the columns of the matrix 
\[
\begin{bmatrix*}[r]
 -1&-4&-3&-2&-6&-7& 0\\
  6&5&1&8&0&0&3
\end{bmatrix*}
\]
and $\mathbb P_\Delta$ is the corresponding
toric surface.
The linear system of Laurent polynomials
whose monomial exponents are in $\Delta$ 
and have multiplicity $7$ at $(1,1)$ 
contains a unique irreducible curve. 
The strict transform of the closure of this
curve is a smooth elliptic curve $C$ of 
$X :=\Bl_e\mathbb P_\Delta$.
It has the minimal equation
\[
 y^2 + xy = x^3 - x^2 - 4x + 4.
\]
This is the curve labeled 
\href{https://www.lmfdb.org/EllipticCurve/Q/446/a/1}{446.a1}
in the LMFDB database
and it has Mordell-Weil group $\mathbb Z^2$.
It is possible to show that 
$\mathcal O_C(C)$ is a non-trivial
element of the Mordell-Weil group
and so it is not torsion
because the Mordell-Weil group is 
torsion-free.
As a consequence $h^0(X,nC) = 1$ for
any $n>0$ so that $[C]$ spans an extremal
ray of $\overline\Eff(X)$. By the 
Riemann-Roch theorem $\overline{\rm Eff}(X)$
must contain the circular cone generated
by the classes of divisors $D$ with $D^2\geq 0$ and $D\cdot H\geq 0$,
where $H$ is an ample class.

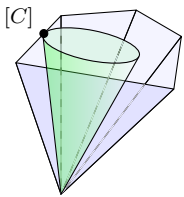
\begin{figure}[ht]
\begin{center}

\begin{tikzpicture}[join=round, scale=.5]
    \tikzstyle{conefill} = [fill=blue!20,fill opacity=0.2]
    \tikzstyle{ann} = [fill=white,font=\footnotesize,inner sep=1pt]
    \tikzstyle{ghostfill} = [fill=white]
         \tikzstyle{ghostdraw} = [draw=black!50]

    \filldraw[conefill](0,0)--(3,2.8)--(3.2,4.2)--cycle;

    \filldraw[conefill](0,0)--(3.2,4.2)--(2,4.9)--cycle;

    \filldraw[conefill](0,0)--(2,4.9)--(0,4.75)--cycle;

    \filldraw[conefill](0,0)--(0,4.75)--(-1.2,3.5)--cycle;

    \filldraw[fill=blue!20,fill opacity=0.45](0,0)--(-1.2,3.5)--(3,2.8)--cycle;

\draw[black!30!white] (0,0)--(0,4.43);   %This is to make backlines lighter

\draw[blue!8!white,dashed] (0,0)--(0,4.43);   %Dashed to the circle

\draw[black!30!white] (0,0)--(2/1.2,4.9/1.2);   %This is to make backlines lighter

\draw[blue!8!white,dashed] (0,0)--(2/1.2,4.9/1.2);    %Dashed to the circle

\draw[black!30!white] (0,0)--(3.2/1.435,4.2/1.435);   %This is to make backlines lighter

\draw[blue!8!white,dashed] (0,0)--(3.2/1.435,4.2/1.435);   %Dashed to the circle

  \def\x{1.3}        % This one changes the radius
  \def\y{4.0}
  \def\R{\x+0.005}
  \def\yc{\y+0.04}
  \def\e{0.4}

  % AK4 cone 1
  \begin{scope}[rotate=-11.5]   % Good for radius   \def\x{1}    
    \shade[right color=white,left color=green,opacity=0.3]
      (-\x,\yc) -- (-\x,\yc) arc (180:360:{\R} and \e) -- (\x,\yc) -- (0,0) -- cycle;
    \draw[fill=green,opacity=0.08]
      (0,\yc) circle ({\R} and \e);
    \draw
      (-\x,\y) -- (0,0) -- (\x,\y);
    \draw
      (0,\yc) circle ({\R} and \e);
  \tkzDefPoint(-\x,\y+0.1){P}
  \tkzDrawPoints[fill=black,color=black,size=3](P)
  \node[above left] at (-\x,\y) {\footnotesize $[C]$}; 
  \end{scope}
\end{tikzpicture}
\end{center}
\caption{A circular cone inside the pseudoeffective cone $\overline{\rm Eff}(X)$. If the class $[C]$ spans an extremal ray of $\overline{\rm Eff}(X)$, then this cone cannot be polyhedral.}
\end{figure}

By a convexity argument one concludes
that $\overline{\rm Eff}(X)$ cannot
be polyhedral.
By Remarks~\ref{remark.effective.cone} 
and \ref{remark.fg.cone}
one concludes that
the effective cone of $\overline M_{0,10}$
is not polyhedral and thus, by the latter
remark, the same holds for
$\overline M_{0,n}$ for $n\geq 10$.

\end{document}